\documentclass[12pt]{amsart}
\pdfoutput=1
\usepackage[margin=1.25in]{geometry}
\usepackage{amsmath}
\usepackage{graphicx}
\usepackage{latexsym}
\usepackage{color}
\usepackage{amscd}
\usepackage[all]{xy}
\usepackage{enumerate}
\usepackage{hyperref}
\usepackage{soul}
\usepackage{comment}
\usepackage{epsfig}
\usepackage{epstopdf}
\usepackage{tikz}
\usepackage{marginnote}
\usepackage{lipsum}
\usepackage{mathtools}
\usepackage{caption} 
\usepackage{subcaption}

\usepackage{multirow}

\newtheorem{thm}{Theorem}

\newtheorem{lemma}[thm]{Lemma}

\newtheorem{proposition}[thm]{Proposition}

\newtheorem*{mthma}{Theorem~A}
\newtheorem*{mthmb}{Theorem~B}

\theoremstyle{definition}

\newtheorem*{definition*}{Definition}

\newcommand{\CPb}{\overline{\mathbb{CP}}{}^{2}}
\newcommand{\CP}{{\mathbb{CP}}{}^{2}}

\newcommand{\R}{\mathbb{R}}

\newcommand{\N}{\mathbb{N}}
\newcommand{\Z}{\mathbb{Z}}

\newcommand{\M}{\operatorname{Mod}}

\def \x {\times}
\def \eu{{\rm{e}}}

 \def\R{{\mathbb{R}}}
 \def\Z{{\mathbb{Z}}}
 
 \def\N{{\mathbb{N}}}
 \def\S{{\Sigma}}

\begin{document}

\title[Spin Lefschetz fibrations are abundant]{Spin Lefschetz fibrations are abundant}

\author[M. Arabadji]{Mihail Arabadji}
\address{Department of Mathematics and Statistics, University of Massachusetts, Amherst, MA 01003, USA}
\email{marabadji@umass.edu}

\author[R. \.{I}. Baykur]{R. \.{I}nan\c{c} Baykur}
\address{Department of Mathematics and Statistics, University of Massachusetts, Amherst, MA 01003, USA}
\email{inanc.baykur@umass.edu}

\maketitle 

\begin{abstract}
We prove that any finitely presented group can be realized as the fundamental group of a spin Lefschetz fibration over the $2$--sphere.  We moreover show that any  admissible lattice point in the symplectic geography plane below the Noether line can be realized by a simply-connected spin Lefschetz fibration.
\end{abstract}

\vspace{0.3in}
\section{Introduction}	

Explicit constructions of Lefschetz fibrations with prescribed fundamental groups were given by Amoros, Bogomolov,  Katzarkov and Pantev \cite{AmorosEtal} and by Korkmaz \cite{Korkmaz}; also see \cite{HamadaEtal}.  We show that the same result holds for a much smaller family of Lefschetz fibrations:

\begin{mthma}
Given any finitely presented group $G$, there exists a spin symplectic Lefschetz fibration $X\to S^2$ with $\pi_1(X) \cong G$. 
\label{thm:1}	
\end{mthma}

These results were inspired by the pioneering work of  Gompf,  who proved that any finitely presented group $G$  is the fundamental group of a closed  symplectic $4$--manifold \cite{Gompf}, which can  moreover be assumed  to be spin.  By the existence of Lefschetz pencils on any symplectic \mbox{$4$--manifold} due to Donaldson \cite{Donaldson},  it then follows a priori that,  after blowing up the base points of the pencil,   one can realize $G$ as the fundamental group of a symplectic Lefschetz fibration, however, these are never spin. 

On the other hand,  unlike K\"{a}hler surfaces,  there are minimal symplectic $4$--manifolds of general type violating the Noether inequality,  which was shown again by Gompf \cite{Gompf}.  More recently,  Korkmaz, Simone, and the second author showed that all the lattice points in the symplectic geography plane below the Noether line  can be further realized by simply-connected symplectic Lefschetz fibrations  \cite{BaykurEtal}.  We prove that a similar result holds in the spin case: 

\begin{mthmb}
For any pair of non-negative integers $(m,n)$ satisfying the inequalities $n\geq 0$,  $n\equiv 8m $ (mod 16),   $n\leq 8(m - 6)$ and $n\leq \frac{16}{3} \, m $,  there exists a simply-connected spin symplectic Lefschetz fibration $X\to S^2$ such that $\chi_h(X)= m$ and $c_1^2(X)=n$.  In particular,  any admissible point in the symplectic geography plane below the Noether line is realized by a simply-connected spin Lefschetz fibration.
\label{thm:2}
\end{mthmb}

Among the hypotheses in the theorem,  the first inequality is due to a theorem of Taubes, who showed that $c_1^2(X) \geq 0$ for any non-ruled minimal symplectic $4$--manifold  $X$,  whereas the second  equality follows from Rokhlin's theorem.  A pair $(m,n) \in \N^2$ satisfying this condition is called \emph{admissible}.  The first systematic production of  spin symplectic $4$--manifolds realizing the above admissible lattice points,  but without the Lefschetz fibration structure we get,  were first obtained by J. Park in \cite{Park}.

Our examples are produced explicitly via positive Dehn twist factorizations in the mapping class group. The spin Lefschetz fibrations for Theorems~A and B are obtained by adapting the strategies of \cite{Korkmaz} and \cite{BaykurEtal}, respectively, together with a subtle use of the breeding technique \cite{BaykurGenus3, BaykurHamada} for the latter.
The main challenge in producing the examples in either theorem is due to the fact that the monodromy of a spin Lefschetz fibration lies in a proper subgroup of the mapping class group (fixing a spin structure on the fiber),  so throughout our work,  we restrain ourselves to algebraic manipulations in this smaller mapping class group. 

\vspace{0.1in}
\noindent \textit{Acknowledgements.} We would like to thank Noriyuki Hamada and Mustafa Korkmaz for their helpful comments on a draft of our paper.  This work was  supported by the NSF grant  DMS-2005327.  The second author would like to thank Harvard University and  Max Planck Institute for Mathematics in Bonn for their hospitality during the writing of this article. 

\section{Preliminaries}  \label{sec:preliminaries}

We begin with a crash review of the concepts and  background results underlying the rest of our article, along with our conventions. We refer the reader to  \cite{GompfStipsicz} more details and comprehensive references on Lefschetz fibrations, symplectic $4$--manifolds, and monodromy factorizations, and to \cite{BaykurHamada} for their interplay with spin structures.

\subsection{Lefschetz fibrations and positive factorizations} \

A \textit{Lefschetz fibration} on a closed smooth oriented $4$--manifold $X$ is a smooth surjective map $f\colon X \to S^2$, a submersion on the complement of finitely many points $\{p_i\}\neq \emptyset$ all in distinct fibers, around which $f$ conforms (compatibly with fixed global orientations on $X$ and $S^2$) to the local complex model of a nodal singularity $(z_1,z_2) \mapsto z_1 z_2$. We assume that there are no exceptional spheres contained in the fibers. Each nodal fiber of the Lefschetz fibration $(X,f)$ is obtained by crashing a simple closed curve,  called a \textit{vanishing cycle},  on a reference regular fiber $F$.

We denote by $\Sigma_g^b$ a compact connected oriented surface of genus $g$ with $b$ boundary components. Let $\textrm{Diff}^+(\Sigma_g^b)$ denote the group of orientation-preserving diffeomorphisms of $\Sigma_g^b$ compactly supported away from the boundary. The \textit{mapping class group} of $\Sigma_g^b$ is defined as $\M(\Sigma_g^b):=\pi_0(\textrm{Diff}^+(\Sigma_g^b))$. When $b=0$, we simply drop $b$ from the above notation. 
Unless mentioned otherwise, by a \emph{curve} $c$ on $\Sigma_g^b$ we mean a smooth simple closed curve.

We denote by $t_c \in \M(\Sigma_g^b)$ the positive (right-handed) Dehn twist along the curve $c \subset \Sigma_g^b$. For any $\psi, \phi  \in \M(\Sigma_g^b)$ we write the conjugate of $\psi$ by $\phi$ as $\psi^\phi=\phi \psi \phi^{-1}$. We act on any curve $c$ in the order $(\varphi \phi)(c)=\varphi(\phi(c))$. An elementary but crucial point is that $t_{c}^{\phi}=t_{\phi(c)}$. For any product of Dehn twists $W=\prod_{i=1}^\ell t_{c_i}^{k_i}$ and $\phi$ in $\M(\Sigma_g^b)$, we denote the conjugated product by $W^\phi = \prod_{i=1}^\ell t_{\phi(c_i)}^{k_i}$.

Let $\{c_i\}$ be a non-empty collection of curves on $\Sigma_g^b$ which do not become null-homotopic after an embedding $\Sigma_g^b  \hookrightarrow\Sigma_g$.  Let $\{\delta_j\}$ be a collection of $b$ curves  parallel to distinct boundary components of $\Sigma_g^b$. A relation of the form
\begin{equation} \label{eqn:PF1}
t_{c_1} t_{c_2} \cdots t_{c_l} = t_{\delta_1}^{k_1} \cdots t_{\delta_b}^{k_b} 
\ \ \ \text{ in } \ \M(\Sigma_g^b)
\end{equation}
corresponds to a genus--$g$ Lefschetz fibration $(X,f)$ with a reference regular fiber $F$ identified with $\Sigma_g$, with vanishing cycles $\{c_i\}$ and $b$ disjoint sections $\{S_j\}$ of self-intersections $S_j \cdot S_j=-k_j$. 

The product on the left-hand side of the equality~\eqref{eqn:PF1}, the word $P$ in positive Dehn twists, is called a \emph{positive factorization} of the mapping class on the right-hand side that maps to the trivial word under the homomorphism induced by an embedding $\Sigma_g^b \hookrightarrow \Sigma_g$. We will often denote the corresponding Lefschetz fibration as $X_P$. 

As shown by Gompf, every Lefschetz fibration $(X,f)$ admits a Thurston type symplectic form with respect to which the fibers are symplectic. 

\subsection{Fiber sums and fundamental groups} \

A  Lefschetz fibration $X_P$ corresponding to a positive factorization $P:=t_{c_1} t_{c_2} \cdots t_{c_l}$ in $\M(\Sigma_g^1)$ (of some power of the boundary twist)  has $\pi_1(X_P) \cong \, \pi_1(\Sigma_g) \, / \, N(\{c_i\})$, where $N(\{c_i\})$ is the subgroup of $\pi_1(\Sigma_g)$  generated normally by collection of the vanishing cycles $c_i$. 

Given $P_1:=t_{c_1} t_{c_2} \cdots t_{c_l} =t_{\delta}^{k_1}$ and $P_2:=
t_{d_1} t_{d_2} \cdots t_{d_l}=t_{\delta}^{k_2}$, and any $\phi \in \M(\Sigma_g^1)$, we can always derive another positive factorization $P_1 P_2^{\phi}= t_{\delta}^{k_1+k_2}$ in $\M(\Sigma_g^1)$, prescribing a new Lefschetz fibration $X_{P_1P_2^\phi}$ with a section of self-intersection $-(k_1+k_2)$. This coincides with the well-known \emph{twisted fiber sum} operation applied to the Lefschetz fibrations $X_{P_1}$ and $X_{P_2}$. We have $\pi_1(X_{P_1P_2^\phi}) \cong \pi_1(\Sigma_g) \, / \, N(\{c_i\} \cup \{\phi(d_j)\})$.  

A neat trick of Korkmaz, applicable in the more special setting described in the next proposition, 
will come very handy for our arguments to follow:

\begin{proposition}[Korkmaz \cite{Korkmaz}]
\label{lem:presentation}
Let $P= t_{c_1}t_{c_2}\cdots t_{c_\ell}$ be a positive factorization of (some power of) a boundary twist in $\M(\Sigma_g^1)$. Let $d$ be a curve on $\Sigma_g$ intersecting at least one $c_i$ transversally at one point. Then $\pi_1(X_{PP^{t_d}}) \cong \pi_1(\Sigma_g)/ N(\{c_i\} \cup \{d\})$.
\end{proposition}

\subsection{Spin monodromies and fibrations} \

A \textit{spin structure} $s$ on $\Sigma_g$ is a cohomology class $s \in H^1(UT(\Sigma_g); \Z_2)$ evaluating to $1$ on a fiber of the unit tangent bundle $UT(\Sigma_g)$. There is a bijection between the set of spin structures on $\Sigma_g$, which we denote by $\textrm{Spin}(\Sigma_g)$, and the set of quadratic forms on $H_1(\Sigma_g;\Z_2)$ with respect to the intersection pairing. Recall that $q\colon H_1(\Sigma_g;\Z_2) \to \Z_2$ is such a quadratic form if $q(a+b)=q(a)+q(b)+a\cdot b$ for every $a, b\in H_1(\Sigma_g; \Z_2)$. 

For a fixed spin structure $s$ on $\Sigma_g$, the \textit{spin mapping class group} $\M(\Sigma_g, s)$ is the stabilizer group of $s$, or equivalently that of the corresponding quadratic form $q$, in $\M(\Sigma_g)$. 
For any non-separating curve $c \subset \Sigma_g$, we have $t_c \in \M(\Sigma_g, s)$ if and only if $q(c)=1$.

The following, which is a reformulation of a theorem of Stipsicz, provides us with a criterion for the existence of a spin structure on a Lefschetz fibration:

\begin{thm}[Stipsicz \cite{Stipsicz}]  \label{SpinLF}
Let $X_P$ be the Lefschetz fibration prescribed by a positive factorization $P:=t_{c_1} t_{c_2} \cdots t_{c_l} =t_{\delta}^{k}$ in $\M(\Sigma_g^1)$, and let us denote the images of the twist curves under the embedding $\Sigma_g^1 \hookrightarrow \Sigma_g$ also by $\{c_i\}$. Then, $X_P$ admits a spin structure with a quadratic form $q$ if and only if $k$ is even and $q(c_i)=1$ for all $i$.
\end{thm}

\section{Spin Lefschetz fibrations with prescribed fundamental group}

In this section, we prove Theorem~A, adapting the strategy in \cite{Korkmaz}, where Korkmaz takes twisted fiber sums of many copies of the same Lefschetz fibration (the building block) to obtain a new Lefschetz fibration whose fundamental group is the prescribed finitely presented group. To accomplish the same  with spin fibrations, there are two essential refinements we will need to make. First is to identify a building block $X_P$ where the monodromy curves in the positive factorization $P$  will satisfy the spin condition for some quadratic form we will describe. That is, we will show that  $P:=t_{c_1} t_{c_2} \cdots t_{c_\ell}$ in $\M(\Sigma_g, s)$ for a carefully chosen spin structure $s$. Second is to make sure that when taking the twisted fiber sums to land on the desired fundamental group, in the corresponding positive factorization $P P^{\phi_1} \cdots P^{\phi_m}$, we only use conjugations $\phi_i \in \M(\Sigma_g, s)$.

\subsection{The building block} \

\begin{figure}[h]
	\begin{tikzpicture}[scale=0.4]
		\begin{scope} 
			
			\draw[very thick,rounded corners=15pt] (-6.5,2.5) --(-7.45, 0)--(-6.5,-2.5) -- (12.5,-2.5)--(12.5,2.5) -- (-5.5,2.5);			
			\draw[very thick, rounded corners = 3pt] (-6.5,2.5)--(-6,2.7)--(-5.4,2.5)--(-6,2.3)--(-6.5,2.5); 

			\draw[very thick, xshift=-5.2cm] (0,0) circle [radius=0.6cm];
			\draw[very thick, xshift=-2.6cm] (0,0) circle [radius=0.6cm];
			\draw[very thick, xshift=2.6cm] (0,0) circle [radius=0.6cm];
			\draw[very thick, xshift=7.8cm] (0,0) circle [radius=0.6cm];
			\draw[very thick, xshift=10.4cm] (0,0) circle [radius=0.6cm];
			\draw[very thick, xshift=-0.2cm] (0,0) circle [radius=0.06cm];
			\draw[very thick, xshift=0.2cm] (0,0) circle [radius=0.06cm];
			\draw[very thick, xshift=0.6cm] (0,0) circle [radius=0.06cm];
			\draw[very thick, xshift=4.6cm] (0,0) circle [radius=0.06cm];
			\draw[very thick, xshift=5cm] (0,0) circle [radius=0.06cm];
			\draw[very thick, xshift=5.4cm] (0,0) circle [radius=0.06cm];
			
			
			\draw[red,rounded corners=8] (-4,2.5)--(-3,1)--(12,1.5)--(12,-1.5) -- (-6,-1.6) --(-7.2,0);
			\draw[red,dashed,rounded corners=10 ] (-7.2,0)--(-6,1.5)--(-4,2.5);
			
			\draw[red,rounded corners=5] (0,-2.5)-- (-5.5,-1.2) --(-5.8,-0.2);
			\draw[red,dashed,rounded corners=10 ] (-5.8,0)--(-6,0.8)--(-1,2.5);
			\draw[red,rounded corners=5] (-1,2.5)--(0,1.7)--(5,1.3)--(11.1,0.9)--(11,0.1);
			\draw[red, dashed,rounded corners = 10] (11,-0.2)-- (11.2,-1)--(2,-2) -- (0,-2.5);
			
			\draw[red ,rounded corners=7] (1.5,-2.5)--(-3.8,-0.7) --(-4.6,0);
			\draw[red ,dashed,rounded corners=10 ] (-4.6,0)--(-1.3,1.4)--(1,2.5);
			\draw[red ,rounded corners=10] (1,2.5)--(5,2)--(9,1)--(9.8,0);
			\draw[red, dashed,rounded corners = 10] (9.8,0) --(8,-1.8)--(4,-2)--(1.5,-2.5);
			
			\draw[red, rounded corners=5] (2.5,-2.5)--(1.9,-1.2)--(2,-0.2);
			\draw[red, rounded corners=5, dashed]  (2,0.2)-- (1.9,1.2) -- (2.3,2.5);
			\draw[red, rounded corners=5] (2.3,2.5) -- (3,1.2)--(3.2,0.2);
			\draw[red, rounded corners=5, dashed] (3.2,-0.2)--(3,-1.2)--(2.5,-2.5);
			
			\draw[red, rounded corners=5] (2.6,0.6)--(2.5,1.6)-- (3,2.5);
			\draw[red, rounded corners = 5, dashed ] (3,2.5) -- (3.2,1.6) -- (2.6,0.6);
			
			\draw[red, rounded corners=5] (2.6,-0.6)--(2.5,-1.4)--(3,-2.5);
			\draw[red, rounded corners=5, dashed] (3,-2.5)--(3.2,-1.4)--(2.6,-0.6);
			
			\node[xshift = -2cm, scale=0.8] at (-1,3.2) {$\delta$};
			\node[xshift = -2.4cm, scale=0.7] at (2,3) {$B_0$};
			\node[xshift = -2.4cm, scale=0.7] at (5,3) {$B_1$};
			\node[xshift = -2.4cm, scale=0.7] at (7,3) {$B_2$};
			\node[xshift = -2.4cm, scale=0.7] at (8.3,3) {$B_g$};
			\node[xshift = -2.4cm, scale=0.7] at (9.3,3) {$a$};
			\node[xshift = -2.4cm, scale=0.7] at (9.3,-3) {$b$};

		\end{scope}
	\end{tikzpicture}
	\caption{The vanishing cycles $B_i, a, b$ of $X_{P_g}$ in $\Sigma_g^1 \subset \Sigma_g$.}
	\label{figure1}
	\end{figure}
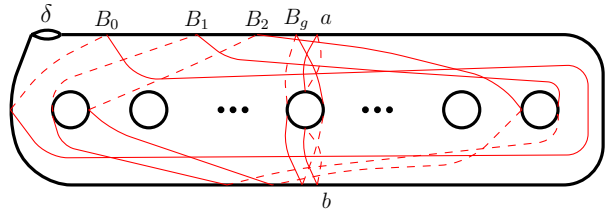

A generalization of the monodromy factorization of the well-known genus--$1$ Lefschetz fibration on $\CP\# 9\, \CPb \cong S^2 \times S^2 \, \# 8\, \CPb$ to any odd genus $g=2n+1$ Lefschetz fibration on $S^2 \times \Sigma_n \, \#8 \, \CPb$ was given by Korkmaz in \cite{Korkmaz0}, and  by Cadavid in \cite{Cadavid}. It has the monodromy factorization:
\begin{equation*}\label{eqn:monodromy}
(t_{B_0}\cdots t_{B_g}t_a^2t_b^2)^2= t_\delta \ \ \ \text{ in } \M(\Sigma_g^1) \, ,
\end{equation*}
where the curves $B_i, a, b$ are shown in the Figure~\ref{figure1}. Capping off the boundary component of $\Sigma_g^1$, we will regard the same curves also in $\Sigma_g$. Let us denote the above positive factorization  by $P_g:=(t_{B_0}\cdots t_{B_g}t_a^2t_b^2)^2$. 

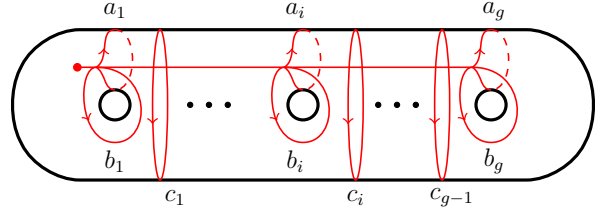
\begin{figure}[h]
	\begin{tikzpicture}[scale=0.5]
		\begin{scope} 
			
			\draw[very thick,rounded corners=15pt] (-7,2) --(-8, 0)--(-7,-2) -- (7,-2)--(8,0) -- (7,2) -- cycle;

			\draw[very thick, xshift=-5cm] (0,0) circle [radius=0.4cm];
			\draw[very thick, xshift=0cm] (0,0) circle [radius=0.4cm];
			\draw[very thick, xshift=5cm] (0,0) circle [radius=0.4cm];
			
			\draw[fill,xshift = -3cm] (0,0) circle [radius = 0.07cm];
			\draw[fill,xshift = -2.5cm] (0,0) circle [radius = 0.07cm];
			\draw[fill,xshift = -2cm] (0,0) circle [radius = 0.07cm];
			\draw[fill,xshift = 3cm] (0,0) circle [radius = 0.07cm];
			\draw[fill,xshift = 2.5cm] (0,0) circle [radius = 0.07cm];
			\draw[fill,xshift = 2cm] (0,0) circle [radius = 0.07cm];
			
			
			\draw[red, fill] (-6,1) circle [radius = 0.1cm];
			\draw[red, semithick] (-6,1)--(4.5,1);
			\draw[red, semithick] (-5.5,1) to [out =30, in = 180] (-5,2);
			\draw[red, semithick, dashed] (-5,2) to [out = 0, in = 0] (-5,0.4);
			\draw[red, semithick] (-5,0.4) to [out=180, in= -30] (-5.5,1);
			\draw[red,semithick, ->] (-5.3,1.5);
			\draw[red, semithick,xshift=5cm] (-5.5,1) to [out =30, in = 180] (-5,2);
			\draw[red, semithick, dashed,xshift=5cm] (-5,2) to [out = 0, in = 0] (-5,0.4);
			\draw[red, semithick,xshift=5cm ] (-5,0.4) to [out=180, in= -30] (-5.5,1);
			\draw[red,semithick, ->, xshift=5cm] (-5.3,1.5);
			\draw[red, semithick,xshift=10cm] (-5.5,1) to [out =30, in = 180] (-5,2);
			\draw[red, semithick, dashed,xshift=10cm] (-5,2) to [out = 0, in = 0] (-5,0.4);
			\draw[red, semithick,xshift=10cm ] (-5,0.4) to [out=180, in= -30] (-5.5,1);
			\draw[red,semithick, ->, xshift=10cm] (-5.3,1.5);
			
			\draw[red, semithick] (-5.5,1) .. controls (-4,1) and (-4,-1) .. (-5,-1) (-5.5,1) .. controls (-6,1) and (-6,-1) .. (-5,-1) ;
			\draw[red, semithick,->] (-5.7,-0.4) --(-5.695,-0.41);
			\draw[red, semithick,xshift = 5cm] (-5.5,1) .. controls (-4,1) and (-4,-1) .. (-5,-1) (-5.5,1) .. controls (-6,1) and (-6,-1) .. (-5,-1) ;
			\draw[red, semithick,->,xshift = 5cm] (-5.7,-0.4) --(-5.695,-0.41);
			\draw[red, semithick, xshift = 10cm] (-5.5,1) .. controls (-4,1) and (-4,-1) .. (-5,-1) (-5.5,1) .. controls (-6,1) and (-6,-1) .. (-5,-1) ;
			\draw[red, semithick,->, xshift = 10cm] (-5.7,-0.4) --(-5.695,-0.41);
			
			\draw[red,semithick  ] (-3.8,0) ellipse (0.2cm and 2cm);
			\draw[red,semithick, ->] (-4,-0.5) -- (-4,-0.51);
			\draw[red,semithick, xshift= 5.2cm  ] (-3.8,0) ellipse (0.2cm and 2cm);
			\draw[red,semithick, ->, xshift = 5.2cm ] (-4,-0.5) -- (-4,-0.51);
			\draw[red,semithick, xshift= 7.5cm  ] (-3.8,0) ellipse (0.2cm and 2cm);
			\draw[red,semithick, ->, xshift = 7.5cm ] (-4,-0.5) -- (-4,-0.51);

			
			\node[scale=0.8] at (-5,2.5) {$a_1$};
			\node[scale=0.8,xshift=3cm] at (-5,2.5) {$a_i$};
			\node[scale=0.8,xshift = 6.3cm ] at (-5,2.5) {$a_g$};
			
			\node[scale=0.8,yshift = -2.5cm] at (-5,2.5) {$b_1$};
			\node[scale=0.8,yshift = -2.5cm, xshift = 3cm] at (-5,2.5) {$b_i$};
			\node[scale=0.8,yshift = -2.5cm, xshift = 6.3cm ] at (-5,2.5) {$b_g$};

			\node[scale=0.8,yshift = -3.1cm,xshift = 1 cm] at (-5,2.5) {$c_1$};
			\node[scale=0.8,yshift = -3.1cm, xshift = 4cm] at (-5,2.5) {$c_i$};
			\node[scale=0.8,yshift = -3.1cm, xshift = 5.6cm ] at (-5,2.5) {$c_{g-1}$};
			
		\end{scope}
	
	\end{tikzpicture}
	\caption{The generators for $\pi_1(\Sigma_g)$ and the $C_i$ curves.}
	\label{figure2}
\end{figure}

Clearly, $\pi_1(X_{P_g}) \cong \pi_1(\Sigma_n)$ will have larger number of generators we can work with as we increase $g=2n+1$. Let us first review the presentation for $\pi_1(X_{P_g})$. Consider the geometric basis $\{a_i, b_i\}_{i=1}^g$ for $\pi_1(\Sigma_g)$, where the based oriented curves $a_i, b_i$ are as shown in Figure~\ref{figure2}. We have 
\[
\pi_1(X) \cong \langle a_1, \dots, a_g,b_1, \ldots, b_g \;|\; C_g, \, a,b, B_0,\dots, B_g  \rangle \, ,
\] 
where\footnote{Here we adopted Korkmaz's generating set to make our calculations comparable to his work in \cite{Korkmaz}, which yields a non-standard expression for the surface relator as iterated conjugates, resulting in $C_g=1$.}
\begin{equation} \label{eqn:pi1expressions}
\begin{cases}	
	B_0 = b_1\cdots b_g \\
	B_{2k-1}= a_k b_k \cdots b_{g+1-k}C_{g+1-k}a_{g+1-k} & 1 \leq k \leq n+1 \\
	B_{2k}= a_k b_{k+1} \cdots b_{g-k} C_{g-k} a_{g+1-k} & 1 \leq k \leq n \\
	a = a_{n+1} \\
	b = C_n a_{n+1} \\
    C_1=b_1^{-1} a_1 b_1 a_1^{-1} \\
    C_i= b_i^{-1}C_{i-1}a_ib_ia_i^{-1}  & 2 \leq i \leq g. 
\end{cases} 
\end{equation}

Next, we will describe a spin structure for which the vanising cycles of this Lefschetz fibration satisfy the monodromy condition.\footnote{It may be worth noting that this is not a trial and error process. By the heuristic arguments of \cite{BaykurHamada}, we are proceeding with an educated guess, since these fibrations are known to come from pencils on the spin manifolds $S^2 \times \Sigma_n$; see \cite{Hamada}.} 
Forgettig the base point, the geometric basis $\{a_j, b_j\}$ for $\pi_1(\Sigma_g)$ in Figure~\ref{figure2} becomes freely homotopic to a standard symplectic basis on $\Sigma_g$. We can then describe a quadratic form with respect to this basis and evaluate it on the mod--$2$ homology classes of the vanishing cycles described in this basis. The latter is easily derived from~\eqref{eqn:pi1expressions}:
\begin{equation*} \label{eqn:H1expressions}
\begin{cases}	
	B_0 = b_1+ \cdots + b_g \\
	B_{2k-1}= a_k  + (b_k + \cdots + b_{g+1-k}) + a_{g+1-k} & 1 \leq k \leq n+1 \\
	B_{2k}= a_k+ (b_{k+1} + \cdots + b_{g-k}) + a_{g+1-k}   & 1 \leq k \leq n \\
	a = a_{n+1} \\
	b = a_{n+1} \, .\\ 
\end{cases} 
\end{equation*}

Set $q\colon H_1(\Sigma_g, \Z_2) \to \Z_2$ as $q(a_i)=q(b_i)=1$. (There are in fact $2^n$ different spin structures that would work here; we are picking the one that will serve our needs the most in the next stages of the proof.) Note that for any ordered set of curves $\{d_j\}$ we have 
$q(\sum_{i=1}^{n} d_i)=\sum_{i=1}^n q(d_i)+\sum_{i<j} 
d_i\cdot d_j$. Thus for each $k$ as above,

\begin{align*}
	q(B_0)=&\;\sum_{i=1}^{g}q(b_i)=g=1\\
	q(B_{2k-1})=&\; q(a_k)+\sum_{i=k}^{g+1-k} q(b_i)+q(a_{g+1-k})+ a_k\cdot b_k + b_{g+1-k}\cdot a_{g+1-k}\\
	=&\;1+(g+1-k-k+1)+1+1+1=1 \\
	q(B_{2k})=&\; q(a_k)+\sum_{i=k+1}^{g-k} q(b_i)+q(a_{g+1-k})\\
	=&\;1+(g-k-k)+1=1   \\
	q(a)= &\; q(a_{n+1})=1 \\
	q(b)=&\; q(a_{n+1})=1.\\
\end{align*} 
Hence all the monodromy curves of $X_{P_g}$ satisfy the spin condition, which is all we needed at this point.\footnote{Recall that $t_\delta$ has odd power in \eqref{eqn:monodromy}, so $X_{P_g}$ is not a spin Lefschetz fibration, as it shouldn't be, remembering that $X_{P_g} \cong S^2 \x \Sigma_n \, \# \, 8 \CPb$.} To sum up, we got

\begin{lemma}\label{lem:legopiece}
Let $s \in \textrm{Spin}(\Sigma_g)$ correspond to the quadratic form with $q(a_i)=q(b_i)=1$, for  $i=1, \ldots, g$, on the symplectic basis $\{a_i, b_i\}$ above. We  have
\begin{equation*}\label{eqn:spinPF}
(t_{B_0}\cdots t_{B_g}t_a^2t_b^2)^2=1 \ \ \ \text{ in } \M(\Sigma_g, s) \, ,
\end{equation*}
where $B_i, a, b$ are the curves on $\Sigma_g^1 \subset \Sigma_g$ in Figure~\ref{figure1}.
\end{lemma}

\subsection{The construction} \

In anticipation of a forthcoming issue, here we deviate a bit from Korkmaz's steps. In order to guarantee that we can represent the relators by embedded curves on $\Sigma_g$, we change the given presentation.
Instead of reinventing the wheel here, we invoke the following result   (cf. \cite{GhigginiEtal}[Lemma~6.2]):
\begin{lemma}[Ghiggini, Golla and Plamanevskaya \cite{GhigginiEtal}] \label{lem:goodpresentation}
	For any finitely presented group $G$, there exists a presentation $G \cong \langle x_1,\dots,x_n \mid r_1, \dots, r_m \rangle$ such that:
	\begin{itemize}
		\item[(i)] each $r_j$ is a positive (no inverses) word in $x_1, \dots, x_n$;
		\item[(ii)] each generator $x_i$ appears at most once in each $r_j$;
		\item[(iii)] the cyclic order (by index) of the generators $x_1, \ldots, x_n$ is preserved in each $r_j$.
	\end{itemize}
\end{lemma}

\noindent This means that for $x_i:=b_i$ in our generating set, we can assume that all the relators in the generating set can be nicely represented by the embedded curves as in Figure~\ref{fig:3} below, where $R_1$ represents $x_2x_3$, $R_2$ represents $x_1x_2 x_4$, and so on. 

\begin{figure}[h]
	\begin{tikzpicture}[scale=0.6]
		
		\draw[very thick,rounded corners=17pt] (7.5,2) --(-7,2) --(-7.5, -0.25)--(-7,-2.5) -- (7.5,-2.5)--(8,-0.25) --cycle;
		
		\draw[very thick, xshift=-5.2cm] (0,0) circle [radius=0.6cm];
		\draw[very thick, xshift=-2.6cm] (0,0) circle [radius=0.6cm];
		\draw[very thick, xshift=0cm] (0,0) circle [radius=0.6cm];
		\draw[very thick, xshift=2.6cm] (0,0) circle [radius=0.6cm];
		\draw[fill] (4.5,0) circle [radius= 0.07];
		\draw[fill] (5,0) circle [radius= 0.07];
		\draw[fill] (5.5,0) circle [radius= 0.07];

		\draw[thick, red, rounded corners=15pt, xshift = 0cm] (-6.7,1)-- (4,1) -- (4,-1.6) --  (1.6,-1.6) -- (1.6,0.8)  -- (-1,0.8) -- (-1,-1.6) --(-6.7,-1.6) -- cycle ;
		\draw[->,very thick,red ] (-3.9,-1.6)--(-3.8,-1.6);

		\draw[thick, blue, rounded corners=5pt, xshift = 0cm] (-3.6,0.8)--(-5.8,0.8) --(-5.8,1.2) -- (1,1.2) -- (1,-1.2) -- (-3.6,-1.2) -- cycle ;
		\draw[->,very thick,blue ] (-1.8,-1.2)--(-1.7,-1.2);
		\draw[fill] (-5.8,1) ellipse (0.1 and 0.1);

		
		\node[scale=1.1] at (0,-2) {$R_1$};
		\node[scale=1.1] at (4,1.5) {$R_2$};
		\node[xshift = -2.4cm, scale=0.9] at (-1.2,-1) {$x_1$};
		\node[scale=0.9] at (-2.6,-0.9) {$x_2$};
		\node[xshift=2cm, scale=0.9] at (-3.3,-0.9) {$x_3$};
		\node[xshift = 4cm, scale=1] at (-4,-1) {$x_4$};

	\end{tikzpicture}
	\caption{Relator curves on $\Sigma_g$. } \label{fig:3}
	
\end{figure}
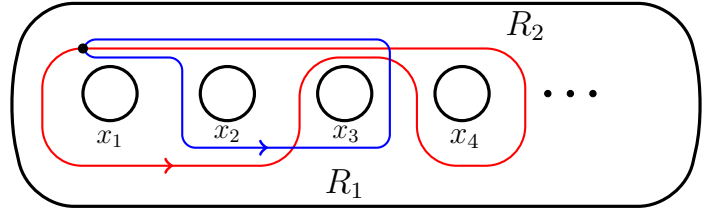

We are now ready to present our construction. 

\begin{proof}[Proof of Theorem~A]
Given a finitely presented group $G$, take a (new) presentation of $G\cong \langle x_1,\dots,x_n \mid r_1, \dots, r_m \rangle$ as in Lemma~\ref{lem:goodpresentation}. Set $g=2n+1$.

Let $P_g:=(t_{B_0}\cdots t_{B_g}t_a^2t_b^2)^2$ be the positive factorization in $\M(\Sigma_g,s)$ given in Lemma~\ref{lem:legopiece}. Because $q(a_i)=1$, we have $t_{a_i} \in \M(\Sigma_g,s)$, for all $i$. So we get a new spin factorization
\begin{equation*}
P_g P_g^{t_{a_1}}P_g^{t_{a_2}} \cdots P_g^{t_{a_g}} =1 \ \ \ \text{ in } \M(\Sigma_g,s) \, 
\end{equation*}
for each odd $g \in \Z^+$, which lifts to a positive factorization of $t_\delta^{g+1}$ in $\M(\Sigma_g^1)$.

From the expression of the monodromy curves of $P_g$ in the $\pi_1(\Sigma_g)$ basis $\{a_i, b_i\}$ given in \eqref{eqn:pi1expressions}, one easily deduces that 
\begin{eqnarray*}
& & 
\langle a_1, \ldots, a_g, b_1, \dots, b_{g} \;|\; C_g \, , a,b,B_0,\dots, B_g, a_1, \ldots, a_g \rangle \, , \\
&\cong& 
\langle b_1, \ldots, b_{2n+1} \;|\; b_1\cdots b_{2n+1}, \;
b_2 \cdots b_{2n},\; \ldots, \; b_n b_{n+1}b_{n+2}, \; b_{n+1} \rangle \, \\
&\cong& 
\langle b_1, \ldots, b_{n}  \rangle \,
\end{eqnarray*}
that is we get a free group on $n$ generators. For the first step, simply note that all $a_j$ and $C_j$ we had in \eqref{eqn:pi1expressions} are trivial in this group.

Now, identifying each generator $x_i$ with $b_i$, for $i=1, \ldots, n$, we can represent each relator $r_j$ by an embedded curve on $R_j$ on $\Sigma_g$. (This is why we switched to this special presentation.) All $\{R_j\}$ can be contained on $\Sigma_n^1 \subset \Sigma_g$ bounded by $c_n$. It is possible that some $q(R_j)=0$. If that is the case, we  replace this $R_j$ with an embedded curve $R'_j$ representing $R_j a_{n+1}$ in $\pi_1(\Sigma_g)$. Such an embedded curve always exists; $R_j$ can be isotoped to meet $a_{n+1}$ only at the base point and one can then resolve the intersection point compatibly with the orientations. So now $q(R'_j)=1$. Otherwise we just take $R'_j:=R_j$. We have $t_{R'_j} \in \M(\Sigma_g,s)$, for all $j=1, \ldots, m$. 

It follows that we have a spin  positive factorization
\begin{equation*}\label{eqn:finalPF}
P_g P_g^{t_{a_1}}P_g^{t_{a_2}} \cdots P_g^{t_{a_g}} \,  P_g^{R'_1} P_g^{R'_2} \cdots P_g^{R'_m} =1 \ \ \ \text{ in } \M(\Sigma_g,s) \, 
\end{equation*}
which now lifts to a positive factorization of $t_\delta^{g+m+1}$ in $\M(\Sigma_g^1)$. If $m$ is odd, we add one more $P_g$ factor to the positive factorization~\eqref{eqn:finalPF}, so then its lift is a positive factorization of $t_\delta^{g+m+2}$. If $m$ is even, leave it as it is. In either case let use denote this final positive factorization in $\M(\Sigma_g,s)$ simply by $P$. Let $X_P$ denote the corresponding Lefschetz fibration.
By Theorem~\ref{SpinLF}, $X_P$ is spin. By Proposition~\ref{lem:presentation}, and the above discussion, we have
\begin{eqnarray*}
\pi_1(X_P)&\cong& 
\langle a_1, \ldots, a_{g},b_1, \ldots, b_{g} \;|\; C_g\, , a,b,B_0,\dots, B_g, a_1, \ldots, a_g, R'_1, \ldots, R'_m \rangle \, , \\
&\cong& 
\langle b_1, b_2, \dots, b_{n} \;|\; R_1, \ldots, R_m \rangle \, \, ,
\end{eqnarray*}
which is the presentation we had for $G$.
\end{proof}

\medskip
\section{Geography of spin Lefschetz fibrations}

In this section we prove Theorem~B by a direct construction of a family of spin Lefschetz fibrations $Z_{g,k}$ populating the region below the Noether line in the geography plane.  We  prescribe these fibrations via new positive factorizations via algebraic manipulations in the mapping class group corresponding to twisted fiber sums and breedings \cite{BaykurGenus3, BaykurHamada}.  We then verify how our careful choice of building blocks out of monodromy factorizations for Lefschetz pencils and fibrations indeed yield positive factorizations in spin mapping class groups.  A somewhat longer calculation will show that our choices also guarantee that $Z_{g,k}$ are simply-connected.  We will then conclude by describing the portion of the geography plane spanned by our spin fibrations.

While some of the particular choices we will make in the construction of  $Z_{g,k}$ may look arbitrary at first,  they are to achieve two somewhat competing properties simultaneously: the existence of a spin structure on $Z_{g,k}$ and the simple-connectivity of $Z_{g,k}$.  The latter calculation implies that the spin structure we describe on $Z_{g,k}$ is in fact unique.

\subsection{The construction} \

\begin{figure}[h!]
	\begin{tikzpicture}[scale=1]
		\begin{scope} [xshift=0cm, yshift=0cm, scale=0.8]
			
			\draw[very thick,rounded corners=22pt] (2.65,1.8) --(-6.5,1.8) --(-7.45, 0)--(-6.5,-1.8) -- (2.65,-1.8)--(3.1,0) --cycle;
			
			\draw[very thick, xshift=-5.2cm] (0,0) circle [radius=0.6cm];
			\draw[very thick, xshift=-2.6cm] (0,0) circle [radius=0.6cm];
			\draw[very thick, xshift=1cm] (0,0) circle [radius=0.6cm];
			
			\filldraw[fill, xshift=-1.2cm] (0,0) circle [radius=0.05cm];
			\filldraw[fill, xshift=-0.8cm] (0,0) circle [radius=0.05cm];
			\filldraw[fill, xshift=-0.4cm] (0,0) circle [radius=0.05cm];
			\draw[thick, red, xshift=-5.2cm] (0,0) circle [radius=0.8cm];
			\draw[thick, red, xshift=-2.6cm] (0,0) circle [radius=0.8cm];
			\draw[thick, red, xshift=1cm] (0,0) circle [radius=0.8cm];
			\draw[thick, red, rounded corners=8pt, xshift=-7.8cm, yshift=0.1cm] (0.6,0) -- (1,-0.13)--(1.7,-0.13) --(2,-0.03);
			\draw[thick,  red, dashed, rounded corners=8pt, xshift=-7.8cm, yshift=0.1cm] (0.75,0.07) -- (0.9,0.13)--(1.7,0.13) --(2,0.03);
			\draw[thick, red, rounded corners=8pt, xshift=-5.2cm, yshift=0.1cm]             (0.6,-0.03) -- (0.9,-0.13)--(1.7,-0.13) --(2,-0.03);
			\draw[thick,  red, dashed, rounded corners=8pt, xshift=-5.2cm, yshift=0.1cm] (0.6,0.03) -- (0.9,0.13)--(1.7,0.13) --(2,0.03);
			\draw[thick, red, rounded corners=8pt, xshift=1cm, yshift=0.1cm]             (0.6,-0.03) -- (0.9,-0.13)--(1.6,-0.13) --(1.9,0);
			\draw[thick,  red, dashed, rounded corners=8pt, xshift=1cm, yshift=0.1cm] (0.6,0.03) -- (0.9,0.13)--(1.7,0.13) --(2,0.03);
			\node[scale=0.8] at (-6.5,-0.3) {$c_1$};
			\node[scale=0.8] at (-3.9,-0.3) {$c_3$};
			\node[scale=0.8] at (2.3,-0.3) {$c_{2g+1}$};
			\node[scale=0.8] at (-5.1,-1.1) {$c_2$};
			\node[scale=0.8] at (-2.5,-1.1) {$c_4$};
			\node[scale=0.8] at (1,-1.1) {$c_{2g}$};
		\end{scope}
	\end{tikzpicture}
	\caption{Dehn twist curves $c_i$ on $\Sigma_g$. } \label{fig:7}
\end{figure}
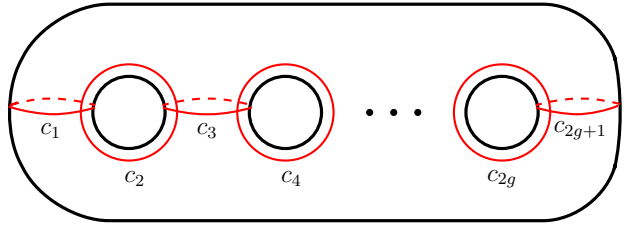

Our first building block is a positive factorization for a Lefschetz fibration on $\CP \# (4g+5) \CPb$ given in \cite{BaykurEtal}. Taking $p=q=2g+2$ in Lemma~4 of  \cite{BaykurEtal},  we obtain
\begin{equation} \label{eqn:hypPF1}
U:= t_1^{2g+2}t_3^{2g+2}(t_1^{t_2}t_2^{t_3}\cdots t_{2g}^{t_{2g+1}})(t_{2g+1}^{t_{2g}}\cdots t_4^{t_3})(t_3^{t_3^{2g+2}t_2}t_2^{t_3^{2g+2}t_1}) =1 
\ \ \ \text{ in }  \M(\S_g)
\end{equation}
which is in fact Hurwitz equivalent to the  square of the positive factorization of the  hyperelliptic involution $h:=t_1\cdots t_{2g}t_{2g+1}^2t_{2g}\cdots t_1$ in $\M(\S_g)$. Here $t_i$ denotes a Dehn twist along the curve $c_i$ shown in Figure~\ref{fig:7}.  We moreover assume that $g \geq 5$ and is odd.  Let also $V$ be the following conjugate of $U$ 
\begin{equation} \label{eqn:hypPF2}
V:= (t_1^{t_2}t_2^{t_3}\cdots t_{2g}^{t_{2g+1}})(t_{2g+1}^{t_{2g}}\cdots t_4^{t_3})(t_3^{t_3^{2g+2}t_2}t_2^{t_3^{2g+2}t_1}) t_1^{2g+2}t_3^{2g+2}=1  \ \ \ \text{ in }  \M(\S_g).
\end{equation}

Consider the following two mapping classes with curves $a,d$ shown in Figure ~\ref{fig:9}:
\begin{eqnarray*}
\phi &:=& (t_8 t_7 t_6 t_a) (t_5 t_6 t_7 t_8) (t_4 t_5 t_6 t_7) (t_3t_4 t_5 t_6) (t_2 t_3 t_4 t_5) (t_1t_2t_3 t_4) \\
\psi &:=& (t_8 t_9 t_{10} t_d) (t_7 t_8 t_9 t_{10})(t_6 t_7 t_8 t_9) \cdots (t_1 t_2 t_3 t_4).
\end{eqnarray*}
We claim that $\phi(c_1)=a$, $\phi(c_3)=b$ and $\psi(c_1)=c$, $\psi (c_3) = d$; see Figure~\ref{fig:9}.  This can be easily verified because of the following elementary observation: Whenever we have a $k$--chain of curves $u_1, \dots u_k$, 
$$t_{u_1}t_{u_2}\cdots t_{u_k}(u_i)=u_{i+1} \;\; \text{ for every } 1\leq i \leq k-1 .$$

Let us denote by $Z_{g}$  the Lefschetz fibration corresponding to the positive factorization $P:=V^\phi U^\psi$ in $\M(\Sigma_g)$, a twisted fiber sum of the Lefschetz fibration on $\CP \# (4g+4) \CPb$ with itself.  Note that we have
\begin{equation} \label{eqn:block1}
	P=V^\phi U^\psi = V_1 \, t_a^{2g+2}t_b^{2g+2} \cdot t_c^{2g+2}t_d^{2g+2} \, U_1 =V_1 \, (t_a t_b t_c t_d)^{2g+2} U_1=1  \ \ \ \text{ in }  \M(\S_g)
\end{equation}
where $V_1, U_1$ are the products of positive Dehn twists
\begin{eqnarray*}
U_1 &:=& ((t_1^{t_2}t_2^{t_3}\cdots t_{2g}^{t_{2g+1}})(t_{2g+1}^{t_{2g}}\cdots t_4^{t_3})(t_3^{t_3^{2g+2}t_2}t_2^{t_3^{2g+2}t_1}) )^\psi \\
V_1 &:=& ((t_1^{t_2}t_2^{t_3}\cdots t_{2g}^{t_{2g+1}})(t_{2g+1}^{t_{2g}}\cdots t_4^{t_3})(t_3^{t_3^{2g+2}t_2}t_2^{t_3^{2g+2}t_1}) )^\phi .
\end{eqnarray*}

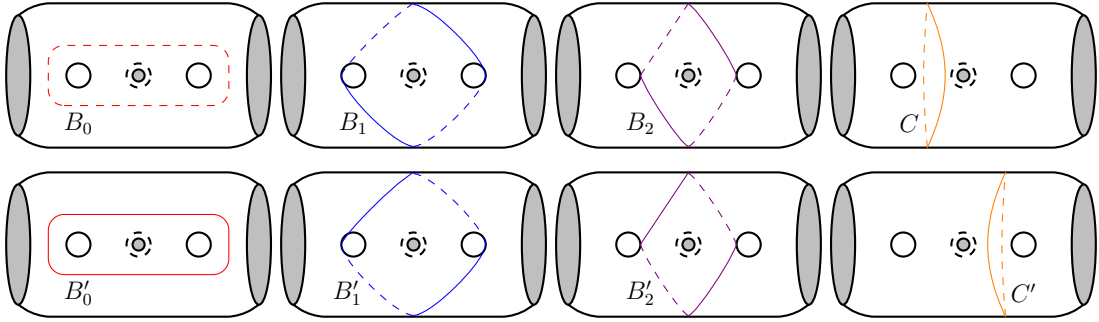
\begin{figure}  
	\begin{tikzpicture}[scale= 0.8]
		\begin{scope}
			
			\draw[thick, rounded corners=5] (0,2)--(0.3,2.2)--(3.7,2.2)--(4,2);
			\draw[thick, rounded corners=5] (0,0)--(0.3,-0.2)--(3.7,-0.2)--(4,0);
			\draw[thick, fill= gray!50] (0,1) ellipse (0.2cm and 1cm);
			\draw[thick, fill = gray!50] (4,1) ellipse (0.2cm and 1cm);
			
			
			\draw[thick] (1,1) ellipse (0.2 and 0.2);
			\draw[thick] (3,1) ellipse (0.2 and 0.2);
			\draw[thick, dashed ] (2,1) ellipse (0.2 and 0.2);
			\draw[thick, fill = gray!50] (2,1) ellipse(0.1 and 0.1);
			
			\draw[red,dashed, rounded corners=6] (0.5,0.5)-- (3.5,0.5)--(3.5,1.5)--(0.5,1.5)--cycle;
			
			\node[scale=0.8] at (1,0.2) {$B_0$};
		\end{scope}	
	\begin{scope}[xshift = 130]
		
		\draw[thick, rounded corners=5] (0,2)--(0.3,2.2)--(3.7,2.2)--(4,2);
		\draw[thick, rounded corners=5] (0,0)--(0.3,-0.2)--(3.7,-0.2)--(4,0);
		\draw[thick, fill= gray!50] (0,1) ellipse (0.2cm and 1cm);
		\draw[thick, fill = gray!50] (4,1) ellipse (0.2cm and 1cm);
		
		
		\draw[thick] (1,1) ellipse (0.2 and 0.2);
		\draw[thick] (3,1) ellipse (0.2 and 0.2);
		\draw[thick, dashed ] (2,1) ellipse (0.2 and 0.2);
		\draw[thick, fill = gray!50] (2,1) ellipse(0.1 and 0.1);
		
		\draw[blue] plot [smooth, tension=1.5] coordinates { (2,-0.18) (1.3,0.3) (0.8,1)};
		\draw[blue,dashed] plot [smooth, tension=1.5] coordinates { (0.8,1) (1.3,1.6)  (2,2.2)};
		\draw[blue] plot [smooth, tension=1.5] coordinates { (2,2.18) (2.7,1.7) (3.2,1)};
		\draw[blue,dashed] plot [smooth, tension=1.5] coordinates { (2,-0.18) (2.7,0.3) (3.2,1)};
		
		\node[scale=0.8] at (1,0.2) {$B_1$};
	\end{scope}

	\begin{scope}[xshift = 260]
		
		\draw[thick, rounded corners=5] (0,2)--(0.3,2.2)--(3.7,2.2)--(4,2);
		\draw[thick, rounded corners=5] (0,0)--(0.3,-0.2)--(3.7,-0.2)--(4,0);
		\draw[thick, fill= gray!50] (0,1) ellipse (0.2cm and 1cm);
		\draw[thick, fill = gray!50] (4,1) ellipse (0.2cm and 1cm);
		
		
		\draw[thick] (1,1) ellipse (0.2 and 0.2);
		\draw[thick] (3,1) ellipse (0.2 and 0.2);
		\draw[thick, dashed ] (2,1) ellipse (0.2 and 0.2);
		\draw[thick, fill = gray!50] (2,1) ellipse(0.1 and 0.1);
		
		\draw[violet] plot [smooth, tension=1.5] coordinates { (2,-0.18) (1.6,0.3) (1.2,1)};
		\draw[violet,dashed] plot [smooth, tension=1.5] coordinates { (1.2,1) (1.6,1.6)  (2,2.2)};
		\draw[violet] plot [smooth, tension=1.5] coordinates { (2,2.18) (2.4,1.7) (2.8,1)};
		\draw[violet,dashed] plot [smooth, tension=1.5] coordinates { (2,-0.18) (2.4,0.3) (2.8,1)};
		
		\node[scale=0.8] at (1.2,0.2) {$B_2$};
	\end{scope}

	\begin{scope}[xshift = 390]
		
		\draw[thick, rounded corners=5] (0,2)--(0.3,2.2)--(3.7,2.2)--(4,2);
		\draw[thick, rounded corners=5] (0,0)--(0.3,-0.2)--(3.7,-0.2)--(4,0);
		\draw[thick, fill= gray!50] (0,1) ellipse (0.2cm and 1cm);
		\draw[thick, fill = gray!50] (4,1) ellipse (0.2cm and 1cm);
		
		
		\draw[thick] (1,1) ellipse (0.2 and 0.2);
		\draw[thick] (3,1) ellipse (0.2 and 0.2);
		\draw[thick, dashed ] (2,1) ellipse (0.2 and 0.2);
		\draw[thick, fill = gray!50] (2,1) ellipse(0.1 and 0.1);
		
		\draw[orange,dashed ] (1.4,-0.2) to [bend right = -5] (1.4,2.2);
		\draw[orange] (1.4,-0.2) to [bend right = 25] (1.4,2.2);

		\node[scale=0.8] at (1.1,0.2) {$C$};	
	\end{scope}

	\begin{scope}[xshift = 260, yshift = -80]
		
		\draw[thick, rounded corners=5] (0,2)--(0.3,2.2)--(3.7,2.2)--(4,2);
		\draw[thick, rounded corners=5] (0,0)--(0.3,-0.2)--(3.7,-0.2)--(4,0);
		\draw[thick, fill= gray!50] (0,1) ellipse (0.2cm and 1cm);
		\draw[thick, fill = gray!50] (4,1) ellipse (0.2cm and 1cm);
		
		
		\draw[thick] (1,1) ellipse (0.2 and 0.2);
		\draw[thick] (3,1) ellipse (0.2 and 0.2);
		\draw[thick, dashed ] (2,1) ellipse (0.2 and 0.2);
		\draw[thick, fill = gray!50] (2,1) ellipse(0.1 and 0.1);
		
		\draw[violet,dashed] plot [smooth, tension=1.5] coordinates { (2,-0.18) (1.6,0.3) (1.2,1)};
		\draw[violet] plot [smooth, tension=1.5] coordinates { (1.2,1) (1.6,1.6)  (2,2.2)};
		\draw[violet,dashed] plot [smooth, tension=1.5] coordinates { (2,2.18) (2.4,1.7) (2.8,1)};
		\draw[violet] plot [smooth, tension=1.5] coordinates { (2,-0.18) (2.4,0.3) (2.8,1)};

		\node[scale=0.8] at (1.2,0.2) {$B_2'$};
	\end{scope}

	\begin{scope}[xshift = 130, yshift = -80]
		
		\draw[thick, rounded corners=5] (0,2)--(0.3,2.2)--(3.7,2.2)--(4,2);
		\draw[thick, rounded corners=5] (0,0)--(0.3,-0.2)--(3.7,-0.2)--(4,0);
		\draw[thick, fill= gray!50] (0,1) ellipse (0.2cm and 1cm);
		\draw[thick, fill = gray!50] (4,1) ellipse (0.2cm and 1cm);
		
		
		\draw[thick] (1,1) ellipse (0.2 and 0.2);
		\draw[thick] (3,1) ellipse (0.2 and 0.2);
		\draw[thick, dashed ] (2,1) ellipse (0.2 and 0.2);
		\draw[thick, fill = gray!50] (2,1) ellipse(0.1 and 0.1);
		
		\draw[blue,dashed] plot [smooth, tension=1.5] coordinates { (2,-0.18) (1.3,0.3) (0.8,1)};
		\draw[blue] plot [smooth, tension=1.5] coordinates { (0.8,1) (1.3,1.6)  (2,2.2)};
		\draw[blue,dashed ] plot [smooth, tension=1.5] coordinates { (2,2.18) (2.7,1.7) (3.2,1)};
		\draw[blue] plot [smooth, tension=1.5] coordinates { (2,-0.18) (2.7,0.3) (3.2,1)};
		
		\node[scale=0.8] at (0.9,0.2) {$B_1'$};
	\end{scope}

	\begin{scope}[xshift = 0, yshift = -80]
		
		\draw[thick, rounded corners=5] (0,2)--(0.3,2.2)--(3.7,2.2)--(4,2);
		\draw[thick, rounded corners=5] (0,0)--(0.3,-0.2)--(3.7,-0.2)--(4,0);
		\draw[thick, fill= gray!50] (0,1) ellipse (0.2cm and 1cm);
		\draw[thick, fill = gray!50] (4,1) ellipse (0.2cm and 1cm);
		
		
		\draw[thick] (1,1) ellipse (0.2 and 0.2);
		\draw[thick] (3,1) ellipse (0.2 and 0.2);
		\draw[thick, dashed ] (2,1) ellipse (0.2 and 0.2);
		\draw[thick, fill = gray!50] (2,1) ellipse(0.1 and 0.1);
		
		\draw[red, rounded corners=6] (0.5,0.5)-- (3.5,0.5)--(3.5,1.5)--(0.5,1.5)--cycle;
		
		\node[scale=0.8] at (1,0.2) {$B_0'$};
	\end{scope}

	\begin{scope}[xshift = 390, yshift = -80]
		
		\draw[thick, rounded corners=5] (0,2)--(0.3,2.2)--(3.7,2.2)--(4,2);
		\draw[thick, rounded corners=5] (0,0)--(0.3,-0.2)--(3.7,-0.2)--(4,0);
		\draw[thick, fill= gray!50] (0,1) ellipse (0.2cm and 1cm);
		\draw[thick, fill = gray!50] (4,1) ellipse (0.2cm and 1cm);
		
		
		\draw[thick] (1,1) ellipse (0.2 and 0.2);
		\draw[thick] (3,1) ellipse (0.2 and 0.2);
		\draw[thick, dashed ] (2,1) ellipse (0.2 and 0.2);
		\draw[thick, fill = gray!50] (2,1) ellipse(0.1 and 0.1);

		\draw[orange,dashed ] (2.7,-0.2) to [bend right = -5] (2.7,2.2);
		\draw[orange] (2.7,-0.2) to [bend right = -25] (2.7,2.2);

		\node[scale=0.8] at (3,0.2) {$C'$};		
	\end{scope}

	\end{tikzpicture}  
	\caption{Vanishing cycles of the genus-2 pencil. }  \label{fig:8}
\end{figure}

Our second building block is the following positive factorization by Hamada 
\begin{equation} \label{eqn:block2}
Q:=t_{B_0}t_{B_1}t_{B_2}t_C t_{C'}t_{B_2'}t_{B_1'}t_{B_0'}=t_{\delta_1}t_{\delta_2}t_{\delta_3}t_{\delta_4} \ \ \ \text{ in }  \M(\Sigma_2^4) 
\end{equation}
for a genus--$2$ Lefschetz pencil on $S^2 \times T^2$,  where the twist curves are as shown in Figure~\ref{fig:8};  see \cite{Hamada, BaykurHamada}. 

Since the curves $\{a,b,c,d\}$  cobound a subsurface $\Sigma_2^4$ of $\Sigma_g$, we can \textit{breed} (see \cite{BaykurGenus3, BaykurHamada}) the genus--$2$ pencil prescribed by~\eqref{eqn:block2} into the Lefschetz fibration prescribed by~\eqref{eqn:block1} for $k$ times,  for any $k  \leq 2g+2$,  and get a new positive factorization
\begin{equation} \label{eqn:monodZ}
	P_{g,k}:=V_1 \, (t_a t_b t_c t_d)^{2g+2-k} R^k U_1=1  \ \ \ \text{ in }  \M(\S_g)
\end{equation}
where $R$ is the image of the positive factorization $Q$ under the homomorphism induced by a specific embedding $\Sigma_2^4 \hookrightarrow \Sigma_g$ we describe below.  We let $Z_{g,k}$ denote the Lefschetz fibration corresponding to the positive factorization $P_{g,k}$.  

The embedding $\Sigma_2^4 \hookrightarrow \Sigma_g$  is described in Figure~\ref{fig:9}.  Brown curves indicate where the boundary curves $\{\delta_i\}$ of $\Sigma_2^4$ in the positive factorization~\eqref{eqn:block2} are mapped to.   Blue arrows illustrate how we isotope the boundaries of $\Sigma_2^4 \subset \R^3$ before embedding it into $\Sigma_g \subset \R^3$.  Red curves constitute a geometric generating set for $H_1(\Sigma_g; \Z_2)$. Red arcs are the parts of these curves contained in the image of the embedding $\Sigma_2^4  \hookrightarrow \Sigma_g$.

 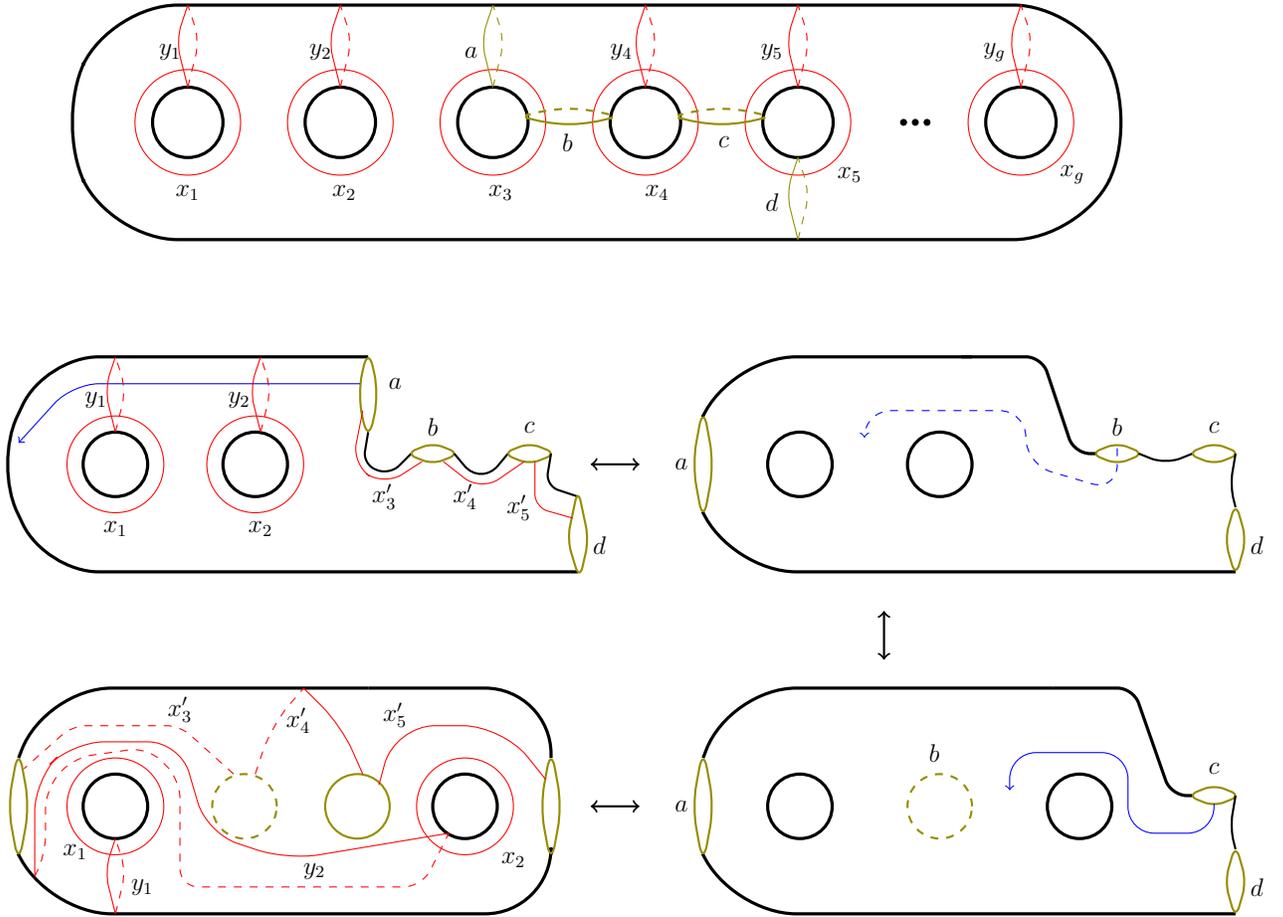
\begin{figure}[h!]
 	\begin{tikzpicture}[scale=1.3]
 		
 		\draw[<->, thick] (2,-1.5)--(2,-2);
 		\draw[<->, thick] (-1,0)--(-0.5,0);
 		\draw[<->, thick] (-1,-3.5)--(-0.5,-3.5);

 		\begin{scope}[xshift = -2cm, yshift=3.5cm, scale=0.6] 
 			
 			
 			\draw[very thick, rounded corners=25pt] (10,2) --(-6.5,2) --(-7.45, 0)--(-6.5,-2) -- (10,-2)--(11,0)--cycle ;

 			\draw[very thick, xshift=-5.2cm] (0,0) circle [radius=0.6cm];
 			\draw[very thick, xshift=-2.6cm] (0,0) circle [radius=0.6cm];
 			\draw[very thick, xshift=0cm] (0,0) circle [radius=0.6cm];
 			\draw[very thick, xshift=+2.6cm] (0,0) circle [radius=0.6cm];
 			\draw[very thick, xshift=+9cm] (0,0) circle [radius=0.6cm];
 			\draw[fill, xshift=7.0cm] (0,0) circle [radius=0.05cm];
 			\draw[fill, xshift=7.2cm] (0,0) circle [radius=0.05cm];
 			\draw[fill, xshift=7.4cm] (0,0) circle [radius=0.05cm];
 			\draw[very thick, xshift=5.2cm] (0,0) circle [radius=0.6cm];
 			
 			\draw[red, xshift=-5.2cm] (0,0) circle [radius=0.9cm] ;
 			\draw[red, xshift=-2.6cm] (0,0) circle [radius=0.9cm] ;
 			\draw[red, xshift=0cm] (0,0) circle [radius=0.9cm] ;
 			\draw[red, xshift=2.6cm] (0,0) circle [radius=0.9cm] ;
 			\draw[red, xshift=5.2cm] (0,0) circle [radius=0.9cm] ;
 			\draw[red, xshift=9cm] (0,0) circle [radius=0.9cm] ;
 			\draw[xshift = -4.7cm, red, rounded corners = 6pt] (-0.5,2)--(-0.7,1.4)--(-0.5,0.6);
 			\draw[xshift = -4.7cm, dashed, red ,rounded corners = 6pt](-0.5,0.6)--(-0.3,1.4)--(-0.5,2);
 			\draw[xshift = -2.1cm, red, rounded corners = 6pt] (-0.5,2)--(-0.7,1.4)--(-0.5,0.6);
 			\draw[xshift = -2.1cm, dashed, red ,rounded corners = 6pt](-0.5,0.6)--(-0.3,1.4)--(-0.5,2);
 			\draw[xshift = 3.1cm, red, rounded corners = 6pt] (-0.5,2)--(-0.7,1.4)--(-0.5,0.6);
 			\draw[xshift = 3.1cm, dashed, red ,rounded corners = 6pt](-0.5,0.6)--(-0.3,1.4)--(-0.5,2);
 			\draw[xshift = 5.7cm, red, rounded corners = 6pt] (-0.5,2)--(-0.7,1.4)--(-0.5,0.6);
 			\draw[xshift = 5.7cm, dashed, red ,rounded corners = 6pt](-0.5,0.6)--(-0.3,1.4)--(-0.5,2);
 			\draw[xshift = 9.5cm, red, rounded corners = 6pt] (-0.5,2)--(-0.7,1.4)--(-0.5,0.6);
 			\draw[xshift = 9.5cm, dashed, red ,rounded corners = 6pt](-0.5,0.6)--(-0.3,1.4)--(-0.5,2);
 			
 			\draw[xshift = 0.5cm, olive, rounded corners = 6pt] (-0.5,2)--(-0.7,1.4)--(-0.5,0.6);
 			\draw[xshift = 0.5cm, dashed, olive ,rounded corners = 6pt](-0.5,0.6)--(-0.3,1.4)--(-0.5,2);
 			\draw[thick, olive, rounded corners=8pt, xshift=0cm, yshift=0.1cm]             (0.6,-0.03) -- (0.9,-0.13)--(1.7,-0.13) --(2,-0.03);
 			\draw[thick,  olive, dashed, rounded corners=8pt, xshift=0cm, yshift=0.1cm] (0.6,0.03) -- (0.9,0.13)--(1.7,0.13) --(2,0.03);
 			\draw[thick, olive, rounded corners=8pt, xshift=2.6cm, yshift=0.1cm]             (0.6,-0.03) -- (0.9,-0.13)--(1.7,-0.13) --(2,-0.03);
 			\draw[thick,  olive, dashed, rounded corners=8pt, xshift=2.6cm, yshift=0.1cm] (0.6,0.03) -- (0.9,0.13)--(1.7,0.13) --(2,0.03);
 			\draw[xshift = 5.7cm, yshift=-2.6cm, olive, rounded corners = 6pt] (-0.5,2)--(-0.7,1.4)--(-0.5,0.6);
 			\draw[xshift = 5.7cm, yshift = -2.6cm, dashed, olive ,rounded corners = 6pt](-0.5,0.6)--(-0.3,1.4)--(-0.5,2);

 			\node[scale=0.8] at (-5.2,-1.2) {$x_1$};
 			\node[scale=0.8] at (-5.5,1.2) {$y_1$};
 			\node[scale=0.8,xshift = 2.6cm] at (-5.2,-1.2) {$x_2$};
 			\node[scale=0.8,xshift = 2.5cm] at (-5.5,1.2) {$y_2$};
 			\node[scale=0.8,,xshift = 5.2cm] at (-5.2,-1.2) {$x_3$};
 			\node[scale=0.8, ,xshift = 7.8cm] at (-5.2,-1.2) {$x_4$};
 			\node[scale=0.8,xshift = 7.5cm] at (-5.5,1.2) {$y_4$};
 			\node[scale=0.8, yshift = 0.3cm,xshift = 11cm] at (-5.2,-1.2) {$x_5$};
 			\node[scale=0.8,xshift = 10cm] at (-5.5,1.2) {$y_5$};
 			\node[scale=0.8, yshift = 0.3cm,xshift = 14.7cm] at (-5.2,-1.2) {$x_g$};
 			\node[scale=0.8,xshift = 13.7cm] at (-5.5,1.2) {$y_g$};

 			\node[scale=0.8,xshift = 5cm] at (-5.5,1.2) {$a$};
 			\node[scale=0.8,xshift = 6.6cm, yshift = -1.5cm] at (-5.5,1.2) {$b$};
 			\node[scale=0.8,xshift = 9.2cm, yshift = -1.5cm] at (-5.5,1.2) {$c$};
 			\node[scale=0.8,xshift = 10cm,yshift = -2.5cm] at (-5.5,1.2) {$d$};

 		\end{scope}

 		\begin{scope} [xshift=-3cm, yshift=0cm, scale=0.55] 

 			\draw[very thick,rounded corners=20pt] (-0.5,2) --(-6.5,2) --(-7.45, 0)--(-6.5,-2) -- (3.4,-2) ;

 			\draw[thick, olive,rounded corners = 8pt] (-0.5,2.2)--(-0.7,1.3)--(-0.5,0.4)--(-0.3,1.3)--cycle; 
 			
 			\draw[thick,rounded corners = 5pt] (-0.5,0.6) -- (-0.6,0)--(-0.1,-0.2)--(0.3,0.2); 
 			
 			\draw[thick, olive, rounded corners = 4pt] (0.2,0.2)-- (0.7,0)--(1.2,0.2)--(0.7,0.4)--cycle; 
 			
 			\draw[thick,rounded corners = 10pt] (1.1,0.2)-- (1.6,-0.4)--(2.1,0.2); 
 			
 			\draw[xshift= 1.8cm, thick, olive, rounded corners = 4pt] (0.2,0.2)-- (0.7,0)--(1.2,0.2)--(0.7,0.4)--cycle; 
 			
 			\draw[thick, rounded corners=5pt] (2.9,0.2)--(2.8,-0.4)--(3.4,-0.6); 
 			
 			\draw[xshift = 3.9cm, yshift = -2.8cm, thick, olive, rounded corners = 5pt] (-0.5,2.35)--(-0.7,1.4)--(-0.47,0.65)--(-0.3,1.4)--cycle; 
 			
 			\draw[very thick, xshift=-5.2cm] (0,0) circle [radius=0.6cm];
 			\draw[very thick, xshift=-2.6cm] (0,0) circle [radius=0.6cm];
 			
 			
 			\draw [->,blue ,rounded corners = 10pt] (-0.65,1.5) -- (-6,1.5)--(-7,0.4);
 			
 			
 			\draw[red, xshift=-5.2cm] (0,0) circle [radius=0.9cm] ; 
 			
 			\draw[red, xshift=-2.6cm] (0,0) circle [radius=0.9cm]; 
 			
 			\draw[red,rounded corners = 8pt] (-0.6,1)--(-0.8,0)--(-0.2,-0.4)--(0.5,0.05); 
 			
 			\draw[red,rounded corners = 8pt] (0.9,0.05)--(1.6,-0.5)--(2.4,0.05); 
 			
 			\draw[red,rounded corners = 5pt] (2.6,0.05)--(2.6,-0.8)--(3.3,-1); 
 			
 			\draw[xshift = -4.7cm, red, rounded corners = 6pt] (-0.5,2)--(-0.7,1.4)--(-0.5,0.6);
 			\draw[xshift = -4.7cm, dashed, red ,rounded corners = 6pt](-0.5,0.6)--(-0.3,1.4)--(-0.5,2); 
 			
 			\draw[xshift = -2cm, red, rounded corners = 6pt] (-0.5,2)--(-0.7,1.4)--(-0.5,0.6);
 			\draw[xshift = -2cm, dashed, red ,rounded corners = 6pt](-0.5,0.6)--(-0.3,1.4)--(-0.5,2); 

 			\node[scale=0.8] at (-5.2,-1.2) {$x_1$};
 			\node[scale=0.8] at (-2.5,-1.2) {$x_2$};
 			\node[scale=0.8] at (-0.2,-0.6) {$x'_3$};
 			\node[scale=0.8] at (1.3,-0.6) {$x'_4$};
 			\node[scale=0.8] at (2.3,-0.8) {$x'_5$};

 			\node[scale=0.8] at (-5.56,1.2) {$y_1$};
 			\node[scale=0.8] at (-2.9,1.2) {$y_2$};

 			\node[scale=0.8] at (0,1.5) {$a$};
 			\node[scale=0.8] at (0.7,0.7) {$b$};
 			\node[scale=0.8] at (2.5,0.7) {$c$};
 			\node[scale=0.8] at (3.8,-1.5) {$d$};

 		\end{scope}

 		\begin{scope} [xshift=4cm, yshift=0cm, scale=0.55]   
 			
 			\draw[very thick,rounded corners=25pt] (-2,2) --(-6.5,2) --(-7, 0.9);
 			\draw[very thick,rounded corners = 6pt] (-2.2,2)--(-0.7,2)--(-0.1,0.2)-- (0.3,0.2);
 			
 			\draw[very thick, rounded corners = 25pt] (-7,-0.9)--(-6.5,-2) -- (2.9,-2)  ;
 			
 			\draw[xshift = -6.5cm, yshift = -1cm,  thick, olive,rounded corners = 8pt] (-0.5,2.1)--(-0.7,1)--(-0.5,-0.1)--(-0.3,1)--cycle;
 			
 			\draw[thick, olive, rounded corners = 4pt] (0.2,0.2)-- (0.7,0)--(1.2,0.2)--(0.7,0.4)--cycle;
 			
 			\draw[thick,rounded corners = 6pt] (1.1,0.2)-- (1.6,0)--(2.1,0.2);
 			
 			\draw[xshift= 1.8cm, thick, olive, rounded corners = 4pt] (0.2,0.2)-- (0.7,0)--(1.2,0.2)--(0.7,0.4)--cycle;
 			
 			\draw[thick, rounded corners=5pt] (2.9,0.2)--(2.8,-0.4)--(2.9,-0.8);
 			
 			\draw[xshift = 3.4cm, yshift = -2.8cm, thick, olive,rounded corners = 5pt] (-0.5,2.1)--(-0.7,1.4)--(-0.5,0.7)--(-0.3,1.4)--cycle;
 			
 			\draw[very thick, xshift=-5.2cm] (0,0) circle [radius=0.6cm];
 			\draw[very thick, xshift=-2.6cm] (0,0) circle [radius=0.6cm];
 			
 			
 			\draw [->,blue,rounded corners = 10pt,dashed ] (0.7,0.3)--(0.7,-0.5)-- (-1,0)--(-1,1)--(-4,1)--(-4,0.5) ;

 			\node[scale=0.8] at (-7.4,0) {$a$};
 			\node[scale=0.8] at (0.7,0.7) {$b$};
 			\node[scale=0.8] at (2.5,0.7) {$c$};
 			\node[scale=0.8] at (3.3,-1.5) {$d$};
 			
 		\end{scope}

 		\begin{scope} [xshift=4cm, yshift=-3.5cm, scale=0.55]   
 			
 			
 			\draw[very thick,rounded corners=25pt] (-0.5,2.2) --(-6.5,2.2) --(-7, 0.9);
 			\draw[very thick,rounded corners = 6pt] (-0.5,2.2)--(1,2.2)--(1.7,0.2)-- (2.1,0.2);
 			\draw[very thick, rounded corners = 25pt] (-7,-0.9)--(-6.5,-2) -- (2.9,-2)  ;
 			
 			\draw[xshift = -6.5cm, yshift = -1cm,  thick, olive,rounded corners = 8pt] (-0.5,2.1)--(-0.7,1)--(-0.5,-0.1)--(-0.3,1)--cycle;
 			
 			\draw[xshift= 1.8cm, thick, olive, rounded corners = 4pt] (0.2,0.2)-- (0.7,0)--(1.2,0.2)--(0.7,0.4)--cycle;
 			
 			\draw[thick, rounded corners=5pt] (2.9,0.2)--(2.8,-0.4)--(2.9,-0.8);
 			
 			\draw[xshift = 3.4cm, yshift = -2.8cm, thick, olive,rounded corners = 5pt] (-0.5,2.1)--(-0.7,1.4)--(-0.5,0.7)--(-0.3,1.4)--cycle;
 			
 			\draw[thick, olive, dashed, xshift=-2.6cm] (0,0) circle [radius=0.6cm];
 			
 			\draw[very thick, xshift=-5.2cm] (0,0) circle [radius=0.6cm];
 			\draw[very thick, xshift=0cm] (0,0) circle [radius=0.6cm];
 			
 			
 			\draw [->,blue ,rounded corners = 10pt] (2.5,0.05)--(2.5,-0.5)--(0.9,-0.5)--(0.9,1)--(-1.3,1)--(-1.3,0.3) ;
 			
 			\node[scale=0.8] at (-7.4,0) {$a$};
 			\node[scale=0.8] at (-2.7,1) {$b$};
 			\node[scale=0.8] at (2.5,0.7) {$c$};
 			\node[scale=0.8] at (3.3,-1.5) {$d$};
 			
 		\end{scope}

 		\begin{scope} [xshift=-3cm, yshift=-3.5cm, scale=0.55]   
 			
 			
 			\draw[very thick,rounded corners=25pt] (-0.5,2.2) --(-6.5,2.2) --(-7, 0.9);
 			\draw[very thick, rounded corners = 25pt] (-7,-0.9)--(-6.5,-2) -- (2.9,-2)--(2.9,-0.9);
 			\draw[very thick,rounded corners = 25pt] (-0.5,2.2)--(2.9,2.2)--(2.9,0.9);
 			
 			\draw[xshift = -6.5cm, yshift = -1cm,  thick, olive,rounded corners = 8pt] (-0.5,2.1)--(-0.7,1)--(-0.5,-0.1)--(-0.3,1)--cycle;

 			\draw[xshift = 3.4cm, yshift = -1cm,  thick, olive,rounded corners = 8pt] (-0.5,2.1)--(-0.7,1)--(-0.5,-0.1)--(-0.3,1)--cycle;
 			
 			\draw[thick, olive, dashed, xshift=-2.8cm] (0,0) circle [radius=0.6cm];
 			
 			\draw[thick, olive, xshift=-0.7cm] (0,0) circle [radius=0.6cm];
 			
 			\draw[very thick, xshift=-5.2cm] (0,0) circle [radius=0.6cm];
 			\draw[very thick, xshift=1.3cm] (0,0) circle [radius=0.6cm];
 			
 			
 			\draw[red, xshift=-5.2cm] (0,0) circle [radius=0.9cm] ; 
 			
 			\draw[red, xshift=1.3cm] (0,0) circle [radius=0.9cm]; 
 			
 			\draw[red, dashed, rounded corners = 5pt] (-3,0.6)-- (-4,1.5)--(-6,1.5)--(-6.9,0.7); 
 			
 			\draw[red, rounded corners = 10pt] (-0.6,0.6)--(-1,1.5)--(-1.7,2.2);
 			\draw[red, dashed, rounded corners = 10pt] (-2.6,0.6)--(-2.3,1.5)--(-1.7,2.2); 
 			
 			\draw[red, rounded corners = 15pt] (-0.3,0.4)--(0,1.5)-- (2,1.5)--(2.8,0.5); 

 			\draw[xshift = -4.7cm, yshift = -2.6cm,  red, rounded corners = 6pt] (-0.5,2)--(-0.7,1.4)--(-0.5,0.6);
 			\draw[xshift = -4.7cm, yshift = -2.6cm, dashed, red ,rounded corners = 6pt](-0.5,0.6)--(-0.3,1.4)--(-0.5,2); 
 			
 			\draw[red, rounded corners = 12pt] (1,-0.5)--(-2,-1)-- (-3.5,-0.5)--(-4,1.2)--(-6,1.2)--(-6.7,0.5)--(-6.7,-1.3);
 			\draw[red, dashed,  rounded corners = 7pt] (-6.7,-1.3)--(-6.5,-0.9)--(-6.5,0.3)--(-6,0.9)--(-4.6,1.1)--(-4,0.5)--(-4,-1.5)--(0.5,-1.5)--(1,-0.5);

 			\node[scale=0.8] at (-5.95,-0.85) {$x_1$};
 			\node[scale=0.8] at (2.2,-1) {$x_2$};
 			\node[scale=0.8] at (-4,1.8) {$x'_3$};
 			\node[scale=0.8] at (-1.8,1.6) {$x_4'$};
 			\node[scale=0.8] at (0,1.7) {$x_5'$};

 			\node[scale=0.8] at (-4.7,-1.5) {$y_1$};
 			\node[scale=0.8] at (-1.5,-1.2) {$y_2$};

 		\end{scope}
 	\end{tikzpicture}
 	\caption{Embedding of $\Sigma_2^4$} \label{fig:9}
 \end{figure}

\subsection{The spin structure on $Z_{g,k}$} \

We are going to invoke Theorem~\ref{SpinLF} to confirm that $Z_{g,k}$ admits a spin structure.  The  curves $\{x_i, y_i\}$ in Figure~\ref{fig:9} constitute a symplectic basis for  $H_1(\Sigma_g ; \Z_2)$.  Consider the following quadratic form $q$ for a spin structure $s \in \textrm{Spin}(\Sigma_g)$ where for any $1 \leq i \leq g$

\begin{equation*}
\begin{cases}
	q(x_i)=1 & \text{ for all } i \\
	q(y_i)=1 & \text{ for } i \text{ odd} \\
	q(y_i)=0 & \text{ for } i \text{ even } \\ 
\end{cases}
\end{equation*}

First of all, $c_{2i}=x_i$, $c_1=y_1$, $c_{2g+1}=y_g$ and $c_{2i+1}=y_i-y_{i+1}$. This means that $q(c_i)=1$ for each $i$.  Therefore,  $t_i:=t_{c_i}  \in \M(\Sigma_g, s)$ for all $i$ and  the positive factorizations given in~\eqref{eqn:hypPF1} and~\eqref{eqn:hypPF2} are in fact factorizations in $\M(\Sigma_g, s)$.  

Secondly,   $a=y_3$ and $d=y_5$ in $H_1(\Sigma_g; \Z_2)$,  so $q(a)=1=q(d)$, in addition to $t_i \in \M(\Sigma_g, s)$, so $\phi,  \psi \in \M(\Sigma_g, s)$.  It follows that $P:=V^\phi U^\psi =1$ is a positive factorization in $\M(\Sigma_g, s)$.

Thirdly,  to check the spin condition for the new monodromy curves in $R$, we would like to express these curves in terms of the generators $\{x_i,y_i\}$.\footnote{Let $F$ denote the embedding $\Sigma_2^4 \hookrightarrow \Sigma_g$ and let $v_j$ be a Dehn twist curve in $R$.  Instead of $F(v_j)\cdot x_i$ and $F(v_j)\cdot y_i$ we can look at $v_j\cdot F^{-1}(x_i)$ and $v_j\cdot F^{-1}(y_i)$ to run the calculation here. Note that if $x_i$ or $y_i$ is only partially contained in the image of $F$, then we denote the arc in its preimage by $x'_i$ or $y'_i$.} We get the following expressions in $H_1(\S_g,\Z_2)$: 
\begin{equation*}
\begin{cases}
	B_0 = x_1+x_2 + y_3+y_4 & B_0'=x_1+x_2+y_4+y_5 \\
	B_1 = x_1+x_2 + y_1+y_2+y_3+y_4+y_5 & B_1' = x_1 + x_2 + y_1 + y_2 + y_4 \\
	B_2 = y_1+y_2+y_3+y_4+y_5 & B_2' = y_1 + y_2 + y_4 \\
	C = y_3 & C' = y_5\\
\end{cases}
\end{equation*}
So we have 
\begin{align*}
	q(B_0)&=q(x_1+x_2 + y_3+y_4) = q(x_1)+q(x_2)+q(y_3)+ q(y_4)=1+1+1+0=1\\
	q(B_0')&=q(x_1+x_2 + y_4+y_5) = q(x_1)+q(x_2)+q(y_4)+ q(y_5)=1+1+0+1=1\\
	q(B_1)&= q(x_1)+q(x_2)+q(y_1)+q(y_2)+q(y_3)+ q(y_4)+q(y_5)+2 \\
	           &=1+1+1+0+ 1+0+1=1 \\
q(B_1')&=q(x_1)+q(x_2)+q(y_1)+q(y_2) + q(y_4) + 2 =1+1+1+0+0=1
\\
	q(B_2)&=q(y_1)+q(y_2)+q(y_3)+ q(y_4)+q(y_5) = 1+0+1+0+1=1\\
	q(B_2')& =q(y_1) +q(y_2)+q(y_4)=1+0+0=1 \\
	q(C)& =q(y_3)=1 \\
	q(C')& =q(y_5)=1.
\end{align*}
Hence, all the vanishing cycles of the Lefschetz fibration $Z_{g,k}$ satisfy the spin condition. 

It is well-known that the Lefschetz fibration with positive factorization $U$ admits a $(-1)$--section; in fact this  fibration is Hurwitz equivalent to a Lefschetz fibration obtained by blowing-up all $4g+4$ base points of a genus--$g$ pencil on $S^2 \x S^2$ \cite{Tanaka}.  Therefore $U,  V$, and in turn $U^\psi, V^\phi$, all lift to a positive factorization of $t_\delta$ in $\M(\Sigma_g^1)$, where $\delta$ is a boundary parallel curve on $\Sigma_g^1$.  We can pick a  $(-1)$--section so that  in the lift of $U$ (and $V$),  the lifts of $t_{c_1},  t_{c_3}$ are still along disjoint curves in $\Sigma_g^1$. \footnote{The Dehn twist curves may get entangled when we take lifts,  but  for just one section we are after, this is not a problem for our positive factorization; see e.g.  \cite{Tanaka} for many possible choices.} The same goes for $t_a,  t_b$ of $V^\phi$ and $t_c, t_d$ of  $U^\psi$.  Let us continue denoting the twist curves in their lifts by $a, b$ and $c,d$.  
After an isotopy,   we can assume that $P=V^\phi U^\psi=1$ lifts to a positive factorization of $t_\delta^2$ in $\M(\Sigma_g^1)$ so that the boundary component is not contained in the subsurface $\Sigma_2^4 \subset \Sigma_g$ cobounded by $\{a,b,c,d\}$.  

Therefore,   for any $k \leq 2g+2$ we have a spin positive factorization 
\begin{equation*} 
P_{g,k}=V_1 \, (t_a t_b t_c t_d)^{2g+2-k} R^k U_1=1  \ \ \ \text{ in }  \M(\S_g, s)
\end{equation*}
which lifts to a positive factorization
 \begin{equation*}
	\tilde{P}_{g,k}=\tilde{V}_1 \, (t_a t_b t_c t_d)^{2g+2-k} \tilde{R}^k \tilde{U}_1=t_\delta^2  \ \ \ \text{ in }  \M(\S_g^1).
\end{equation*}
Hence every $Z_{g,k}$ admits a spin structure by Theorem~\ref{SpinLF}.

\smallskip
\subsection{The fundamental group.} \

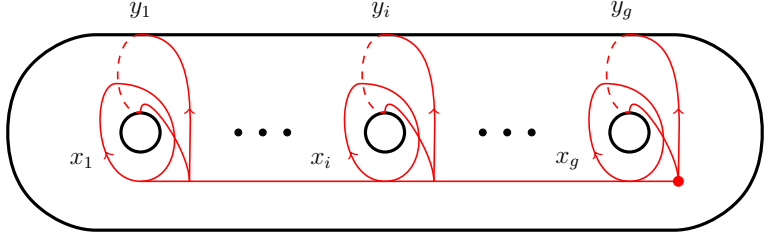
\begin{figure}[h]
	\begin{tikzpicture}[scale=0.65] 
		\begin{scope} 
			
			\draw[very thick,rounded corners=20pt] (-7,2) --(-8, 0)--(-7,-2) -- (7,-2)--(8,0) -- (7,2) -- cycle;

			\draw[very thick, xshift=-5cm] (0,0) circle [radius=0.4cm];
			\draw[very thick, xshift=0cm] (0,0) circle [radius=0.4cm];
			\draw[very thick, xshift=5cm] (0,0) circle [radius=0.4cm];
			
			\draw[fill,xshift = -3cm] (0,0) circle [radius = 0.07cm];
			\draw[fill,xshift = -2.5cm] (0,0) circle [radius = 0.07cm];
			\draw[fill,xshift = -2cm] (0,0) circle [radius = 0.07cm];
			\draw[fill,xshift = 3cm] (0,0) circle [radius = 0.07cm];
			\draw[fill,xshift = 2.5cm] (0,0) circle [radius = 0.07cm];
			\draw[fill,xshift = 2cm] (0,0) circle [radius = 0.07cm];
			
			
			\draw[red, fill] (6,-1) circle [radius = 0.1cm];
			\draw[red, semithick] (-5,-1)--(6,-1);
			
			\draw[red, semithick] (-4,-1) to [out=90, in=90] (-5,0.4);
			\draw[red, semithick] (-4,-1) to [out=90, in=0] (-5,2);
			\draw[red, semithick, dashed] (-5,2) to [out = 180, in = 180] (-5,0.4);
			\draw[red,semithick, ->] (-4,0.5);
			
			\draw[red, semithick, xshift = 5cm] (-4,-1) to [out=90, in=90] (-5,0.4);
			\draw[red, semithick, xshift = 5cm] (-4,-1) to [out=90, in=0] (-5,2);
			\draw[red, semithick, dashed, xshift = 5cm] (-5,2) to [out = 180, in = 180] (-5,0.4);
			\draw[red,semithick, ->, xshift = 5cm] (-4,0.5);
			\draw[red, semithick, xshift = 10cm] (-4,-1) to [out=90, in=90] (-5,0.4);
			\draw[red, semithick, xshift = 10cm] (-4,-1) to [out=90, in=0] (-5,2);
			\draw[red, semithick, dashed, xshift = 10cm] (-5,2) to [out = 180, in = 180] (-5,0.4);
			\draw[red,semithick, ->, xshift = 10cm] (-4,0.5);
			
			\draw[red, semithick] (-5.5,1) .. controls (-4,1) and (-4,-1) .. (-5,-1) (-5.5,1) .. controls (-6,1) and (-6,-1) .. (-5,-1) ;
			\draw[red, semithick,<-] (-5.7,-0.4) --(-5.695,-0.41);
			\draw[red, semithick,xshift = 5cm] (-5.5,1) .. controls (-4,1) and (-4,-1) .. (-5,-1) (-5.5,1) .. controls (-6,1) and (-6,-1) .. (-5,-1) ;
			\draw[red, semithick,<-,xshift = 5cm] (-5.7,-0.4) --(-5.695,-0.41);
			\draw[red, semithick, xshift = 10cm] (-5.5,1) .. controls (-4,1) and (-4,-1) .. (-5,-1) (-5.5,1) .. controls (-6,1) and (-6,-1) .. (-5,-1) ;
			\draw[red, semithick,<-, xshift = 10cm] (-5.7,-0.4) --(-5.695,-0.41);

			\node[scale=0.8] at (-5,2.5) {$y_1$};
			\node[scale=0.8,xshift=4cm] at (-5,2.5) {$y_i$};
			\node[scale=0.8,xshift = 8cm ] at (-5,2.5) {$y_g$};
			
			\node[scale=0.8, yshift = -2.5cm] at (-6.2,2.5) {$x_1$};
			\node[scale=0.8,yshift = -2.5cm, xshift = 3cm] at (-5,2.5) {$x_i$};
			\node[scale=0.8,yshift = -2.5cm, xshift = 7.1cm ] at (-5,2.5) {$x_g$};

		\end{scope}
		
	\end{tikzpicture}
	\caption{The generators $x_i, y_i$ for $\pi_1(\Sigma_g)$.}
	\label{fig:newpi1basis}
\end{figure}

Let $\{ x_i, y_i\}$ be a geometric basis for $\pi_1(\Sigma_g)$ as shown in Figure~\ref{fig:newpi1basis}.  Since $Z_{g,k}$ has a section,  we have $G:=\pi_1(Z_{g,k}) \cong \pi_1(\Sigma_g) \, / \, N(\{ v_j\})$ where $v_j$ are the Dehn twist curves in the  positive factorization $P_{g,k}$ of $Z_{g,k}$.

Set $S:=(t_1^{t_2}t_2^{t_3}\cdots t_{2g}^{t_{2g+1}})(t_{2g+1}^{t_{2g}}\cdots t_4^{t_3})(t_3^{t_3^{2g+2}t_2}t_2^{t_3^{2g+2}t_1})$. So   $U^\psi= t_c^{2g+2}t_d^{2g+2} U_1$ with $U_1=S^\psi$, and $V^\phi=V_1   t_a^{2g+2}t_b^{2g+2}$ with $V_1=S^\phi$.  While the fundamental group of $Z_{g,k}$ can be calculated from  the factorization  $P_{g,k}=S^{\phi}\,(t_at_bt_ct_d)^{2g+2-k}R^kS^{\psi}$,  it can also be calculated from  $P^{\phi^{-1}}= S\,(t_1t_3t_{\phi^{-1}(c)}t_{\phi^{-1}(d)})^{2g+2-k}(R^{\phi^{-1}})^kS^{\phi^{-1}\psi}$.  We will run our calculations for the latter.

For $k < 2g+2$ the Dehn twist curves of the latter factorization contains all the vanishing cycles $\{c_i\}$ in $U$, which we know  kill all the generators of $\pi_1(\Sigma_g)$ to yield trivial the fundamental group,  as $P_U$ has total space $\CP \# (4g+4 ) \CPb$, a simply-connected space.  Thus,  $\pi_1(Z_{g,k})=1$ for any $k < 2g+2$.

For $k=2g+2$, first note that we can connect the vanishing cycles $\{c_i\}$ of $U$ or $V$ to the basepoint (where any two different paths connecting them to the base point will yield the same normal generating set) so that in $\pi_1(\Sigma_g)$ we have $c_{2i}=x_i$, $c_1=y_1$, $c_{2g+1}=y_g$ and $c_{2i+1}=y_i^{-1}x_{i+1}y_{i+1}x_{i+1}^{-1}$ for each $i$.  It follows that $\pi_1(\Sigma_g)$ is generated by $\{c_i\}$.   To get $G$ we quotient $\pi_1(\Sigma_g)$ by normally generated subgroup by relators coming from the Dehn twist curves in $S$ (and not $U$), which are of the form $t_i(c_{i-1})$ with $2\leq i \leq 2g+1$ and $t_3(c_4)$,  along with several other relations.  We may assume that $c_i$ are oriented so that $c_{i-1}\cdot c_{i}= +1$ for all $i$. Then we have the relators $t_i(c_{i-1})=c_{i-1}c_i^{-1}=1$ and $t_3(c_4)=c_4c_3=1$. These relations imply that
$$c_1=c_2=c_3=c_4=\cdots = c_{2g}=c_{2g+1}$$
and 
$$c_3=c_4^{-1}.$$
We thus see that  $G\cong \left<c_1 \; | \; c_1^2, \text{ rest of the relators coming from other vanishing cycles}  \right>$ for our positive factorization  $S(R^{\phi^{-1}})^{2g+2}S^{\phi^{-1}\psi}.$

At this point $G$ is a quotient of the abelian group $\Z_2$ generated by $c_1$,  so it is certainly an abelian group,   and it suffices to show that $H_1(Z_{g,2g+2})=0$.  

We will argue this by observing that the vanishing cycle coming from $t_{\phi^{-1}(B_2)}$ induces a relator killing the homology class of $c_1$.  This is because it is homologous to an odd factor of $c_1$. 
For this reason, it is in fact enough to consider  $\phi^{-1}(B_2)$ in 
$H_1(Z_{g,2g+2}; \Z_2)$.  By the previous computations, we have 
 $$B_2=(c_1)+(c_1+c_3)+\dots + (c_1+c_3+\dots c_9)=c_1+c_5+c_9.$$ Let's apply $\phi^{-1}$. Then 
\begin{align*}
	c_1+c_5+c_9 &\stackrel{t_8^{-1}}{\mapsto} c_1+c_5+c_8+c_9\stackrel{t_7^{-1}}{\mapsto}c_1+c_5+c_7+c_8+c_9\stackrel{t_6^{-1}}{\mapsto} c_1+c_5+c_7+c_8+c_9\\ &\stackrel{t_a^{-1}}{\mapsto} c_1+c_5+c_7+c_8+c_9\stackrel{t_5^{-1}}{\mapsto} c_1+ c_5+ c_7 +c_8 +c_9 \stackrel{t_6^{-1}}{\mapsto} c_1+ c_5+ c_7+ c_8+ c_9\\
	&\stackrel{t_7^{-1}}{\mapsto} c_1+ c_5+ c_8 +c_9 \stackrel{t_8^{-1}}{\mapsto} c_1+ c_5 +c_9\stackrel{t_4^{-1}}{\mapsto} c_1+ c_4 +c_5 +c_9 \stackrel{t_5^{-1}}{\mapsto} c_1 +c_4 +c_9\\ &\stackrel{t_6^{-1}}{\mapsto} c_1+c_4 + c_9\stackrel{t_7^{-1}}{\mapsto} c_1+ c_4 + c_9 \stackrel{t_3^{-1}}{\mapsto} c_1 + c_3 + c_4 + c_9 \stackrel{t_4^{-1}}{\mapsto} c_1 + c_3 + c_9\\ &\stackrel{t_5^{-1}}{\mapsto} c_1 + c_3 + c_9 \stackrel{t_6^{-1}}{\mapsto} c_1 + c_3 + c_9\stackrel{t_2^{-1}}{\mapsto} c_1 + c_3 + c_9 \stackrel{t_3^{-1}}{\mapsto} c_1 + c_3 + c_9\\ &\stackrel{t_4^{-1}}{\mapsto} c_1 + c_3 + c_4 + c_9 \stackrel{t_5^{-1}}{\mapsto} c_1 + c_3 + c_4 + c_5 + c_9
	\stackrel{t_1^{-1}}{\mapsto}c_1 + c_3 + c_4 + c_5 + c_9\\
	&\stackrel{t_2^{-1}}{\mapsto} c_1 + c_3 + c_4 + c_5 + c_9 \stackrel{t_3^{-1}}{\mapsto} c_1 + c_4 + c_5 + c_9 \stackrel{t_4^{-1}}{\mapsto} c_1 + c_5 + c_9 \, \
\end{align*}
which means $\phi^{-1}(B_2)=c_1$ in this abelian group.  Hence $G\cong 1$.

\smallskip
\subsection{The geography} \

We are left with determining the portion of the geography plane populated by our simply-connected spin Lefschetz fibrations 
\[ \{Z_{g,k} \ | \ g\geq 5 \text{ and odd},  k \leq 2g+2 \text{ and non-negative} \}. \]

The Euler characteristic of $Z_{g,k}$ is given by the formula
\[ \eu(Z_{g,k}) = 4 -4g + \ell =  4-4g + (16g+8 +4k)=12(g+1)+4k \, , \]
where $\ell$ is the number of Dehn twists in $P_{g,k}$.

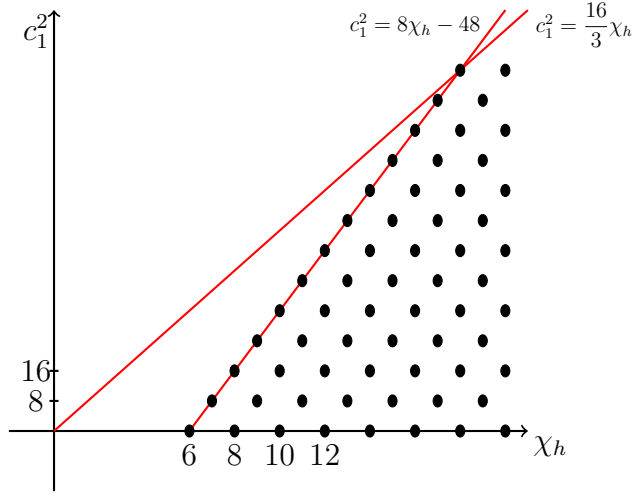
\begin{figure}
	\begin{tikzpicture}[yscale=0.4, xscale=0.3]

		\begin{scope}
			\draw[thick, ->] (-2,0)--(21,0);
			\draw[thick, ->]  (0,-2)--(0,14);
			\draw[thick, color=red] (6,0)--(20,14);
			\draw[thick, red] (0,0)--(21,14);
			
			
			\draw [fill] (6,0) circle [radius=2mm];
			\draw [fill] (8,0) circle [radius=2mm];
			\draw [fill] (10,0) circle [radius=2mm];
			\draw [fill] (12,0) circle [radius=2mm];
			\draw [fill] (14,0) circle [radius=2mm];
			\draw [fill] (16,0) circle [radius=2mm];
			\draw [fill] (18,0) circle [radius=2mm];
			\draw [fill] (20,0) circle [radius=2mm];
			
			\draw [fill] (7,1) circle [radius=2mm];
			\draw [fill] (9,1) circle [radius=2mm];
			\draw [fill] (11,1) circle [radius=2mm];
			\draw [fill] (13,1) circle [radius=2mm];
			\draw [fill] (15,1) circle [radius=2mm];
			\draw [fill] (17,1) circle [radius=2mm];
			\draw [fill] (19,1) circle [radius=2mm];
			
			\draw [fill] (8,2) circle [radius=2mm];
			\draw [fill] (10,2) circle [radius=2mm];
			\draw [fill] (12,2) circle [radius=2mm];
			\draw [fill] (14,2) circle [radius=2mm];
			\draw [fill] (16,2) circle [radius=2mm];
			\draw [fill] (18,2) circle [radius=2mm];
			\draw [fill] (20,2) circle [radius=2mm];
			
			\draw [fill] (9,3) circle [radius=2mm];
			\draw [fill] (11,3) circle [radius=2mm];
			\draw [fill] (13,3) circle [radius=2mm];
			\draw [fill] (15,3) circle [radius=2mm];
			\draw [fill] (17,3) circle [radius=2mm];
			\draw [fill] (19,3) circle [radius=2mm];
			
			\draw [fill] (10,4) circle [radius=2mm];
			\draw [fill] (12,4) circle [radius=2mm];
			\draw [fill] (14,4) circle [radius=2mm];
			\draw [fill] (16,4) circle [radius=2mm];
			\draw [fill] (18,4) circle [radius=2mm];
			\draw [fill] (20,4) circle [radius=2mm];
			
			\draw [fill] (11,5) circle [radius=2mm];
			\draw [fill] (13,5) circle [radius=2mm];
			\draw [fill] (15,5) circle [radius=2mm];
			\draw [fill] (17,5) circle [radius=2mm];
			\draw [fill] (19,5) circle [radius=2mm];
			
			\draw [fill] (12,6) circle [radius=2mm];
			\draw [fill] (14,6) circle [radius=2mm];
			\draw [fill] (16,6) circle [radius=2mm];
			\draw [fill] (18,6) circle [radius=2mm];
			\draw [fill] (20,6) circle [radius=2mm];
			
			\draw [fill] (13,7) circle [radius=2mm];
			\draw [fill] (15,7) circle [radius=2mm];
			\draw [fill] (17,7) circle [radius=2mm];
			\draw [fill] (19,7) circle [radius=2mm];
			
			\draw [fill] (14,8) circle [radius=2mm];
			\draw [fill] (16,8) circle [radius=2mm];
			\draw [fill] (18,8) circle [radius=2mm];
			\draw [fill] (20,8) circle [radius=2mm];
			
			\draw [fill] (15,9) circle [radius=2mm];
			\draw [fill] (17,9) circle [radius=2mm];
			\draw [fill] (19,9) circle [radius=2mm];

			\draw [fill] (16,10) circle [radius=2mm];
			\draw [fill] (18,10) circle [radius=2mm];
			\draw [fill] (20,10) circle [radius=2mm];
			
			\draw [fill] (17,11) circle [radius=2mm];
			\draw [fill] (19,11) circle [radius=2mm];

			\draw [fill] (18,12) circle [radius=2mm];
			\draw [fill] (20,12) circle [radius=2mm];

			\node at (-0.8,13.3) {$c_1^2$};
			\node at (22,-0.5) {$\chi_h$};
			\node at (6,-0.8) {6};
			\node at (8,-0.8) {8};
			\node at (10,-0.8) {10};
			\node at (12,-0.8) {12};
			
			\node at (-0.8,1) {8};
			\draw[thick] (-0.2,1)--(0.2,1);
			\node at (-0.8,2) {16};
			\draw[thick] (-0.2,2)--(0.2,2);
			
			\node[scale=0.7] at (23.5,13.5) {$c_1^2= \dfrac{16}{3}\chi_h$};
			\node[scale=0.7] at (16,13.5) {$c_1^2= 8\chi_h-48$};
			
		\end{scope}
	\end{tikzpicture}
	\caption{The region populated by spin $Z_{g,k}$. } \label{fig:10}
\end{figure}

Since the positive factorization $U$ commutes with a hyperelliptic involution on $\Sigma_g$ (after all, it is Hurwitz equivalent to the positive factorization of a hyperelliptic involution itself),  by Endo's signature formula for hyperelliptic fibrations \cite{Endo},  it has signature $-4g-4$ (as expected, since the total space is $\CP \# (4g+5) \CPb$).  By the Novikov additivity,  we then get $\sigma(Z_g)=-8g-8$.  Breeding the signature zero genus--$2$ Lefschetz pencil into this fibration (any number of times) does not change the signature \cite{BaykurGenus3} and we get  
\[ \sigma(Z_{g,k})= -8(g+1) \,  .\]
We thus have 
\[ 
\chi_h(Z_{g,k})= \frac{1}{4}(\eu(Z_{g,k})+\sigma(Z_{g,k})) = g+1 +k  
\]
and
\[
c_1^2(Z_{g,k}) = 2 \eu(Z_{g,k}) + 3 \sigma(Z_{g,k}) =  8k. 
\]

Thus,  setting $g=2r+5$ with $r\geq 0$,  we see that $\{(\chi_h, c_1^2)(Z_{g,k})\}$ populate the region  
\[ 
\mathcal{R}=\{ (6,0) + r (2,0) + k(1,8) \, | \,  (r, k) \in \N^2 \text{ with } k \leq 4(r+ 3) \}. 
\]

of the geography plane,  or equivalently , 
\[
\mathcal{R} =  \{ (m, n) \in \N^2 \, | \,  n \geq 0 \, , n \leq 8(m-6), \, n \leq \frac{16}{3} m \, \text{ and }  \, n\equiv8m\text{ mod } 16  \}. 
\]
See Figure~\ref{fig:10}.  In particular,  one can easily see from the first description of $\mathcal{R}$ above that we cover all of the admissible lattice points in $\N^2$ under the Noether line.

\end{document}